\documentclass{article}
\usepackage{amsmath,amssymb,amscd}
\usepackage{amsthm}
\usepackage{mathptmx}
\newtheorem{Theorem}{\indent \sc Theorem}[section]
\newtheorem{Proposition}[Theorem]{\indent \sc Proposition}
\newtheorem{Corollary}[Theorem]{\indent \sc Corollary}
\newtheorem{Lemma}[Theorem]{\indent \sc Lemma}
\newtheorem{Remark}[Theorem]{\indent \sc Remark}
\newtheorem{Definition}[Theorem]{\indent \sc Definition}

\newcommand{\spec}{\operatorname{Spec}}
\newcommand{\et}{\mathrm{\acute{e}t}}
\newcommand{\Nis}{\mathrm{Nis}}
\newcommand{\Zar}{\mathrm{Zar}}

\newcommand{\HO}{\operatorname{H}}

\usepackage[all]{xy}

\begin{document}
\author{Makoto Sakagaito}
\title{
On 
\'{e}tale hypercohomology
of henselian regular local rings
with
values in
$p$-adic \'{e}tale Tate twists
}
\date{}
\maketitle
\begin{center}
Indian Institute of Science Education and Research, Bhopal
\\
\textit{E-mail address}: sakagaito43@gmail.com
\end{center}
\begin{abstract}
Let $R$ be the henselization of a local ring of a semistable family
over the spectrum of a discrete valuation ring of mixed characteristic $(0, p)$
and $k$ the residue field of $R$.
In this paper,
we prove an isomorphism of \'{e}tale hypercohomology groups
\begin{math}
\operatorname{H}^{n+1}_{\mathrm{\acute{e}t}}(R, \mathfrak{T}_{r}(n))
\simeq
\operatorname{H}^{1}_{\mathrm{\acute{e}t}}(k, W_{r}\Omega_{\log}^{n})
\end{math}
for any integers $n\geq 0$ and $r>0$ where 
$\mathfrak{T}_{r}(n)$ is the $p$-adic Tate twist and
$W_{r}\Omega_{\log}^{n}$ is the logarithmic Hodge-Witt sheaf. 
As an application, we prove the local-global
principle for Galois cohomology groups over function fields of 
curves over an excellent henselian discrete valuation
ring of mixed characteristic.
\end{abstract}
\section{Introduction}
The main objective of this paper is to study a certain \'{e}tale hypercohomology in the mixed characteristic cases.
In order to understand a scheme $\mathfrak{X}$ over the spectrum of a discrete valuation ring, 
it is effective to observe the special fiber of $\mathfrak{X}$.
So we start in a positive characteristic situation.

Let $k$ be a field of positive characteristic, 
$X$ a \textit{normal crossing variety} over $k$, 
i.e. 
a pure-dimensional scheme of finite type over 
$\spec (k)$ which is everywhere \'{e}tale locally
isomorphic to
\begin{equation*}
\spec (
k[T_{0}, \cdots, T_{d}]
/(T_{0}T_{1}\cdots T_{r})
) 
\end{equation*}
for some integer $r$ with $0\leq r\leq d=\operatorname{dim} X$
(cf.\cite[p.707]{SaL}).

Let us denote
the set of points on $X$
of codimension $i$ by $X^{(i)}$.
For a point $x\in X$, let $i_{x}$ be the canonical map
$x\hookrightarrow X$.
Then
we have three kinds of 
logarithmic Hodge-Witt sheaves:
$W_{r}\omega^{n}_{X, \log}$ (\cite{H}),
\begin{equation*}
\lambda^{n}_{X, r}
:=
\operatorname{Im}\left(
d\log:
(\mathbb{G}_{m, X})^{\otimes n}
\to
\bigoplus_{x\in X^{(0)}}i_{x*}\operatorname{W}_{r}\Omega_{x, \log}^{n}
\right)  
\end{equation*}
and
\begin{equation*}
\nu^{n}_{X, r}:=\operatorname{Ker}
\left(
\partial: 
\bigoplus_{x\in X^{(0)}}
i_{x *}W_{r}\Omega^{n}_{x, \log}
\to
\bigoplus_{x\in X^{(1)}}
i_{x *}W_{r}\Omega^{n-1}_{x, \log}
\right)   
\end{equation*}
(\cite[p.527, 2.2]{SaP}) which agree with
the logarithmic Hodge-Witt sheaves 
$W_{r}\Omega^{n}_{X, \log}$ (\cite{I})
in the case where $X$ is a smooth variety
(cf.\cite[p.528, Th\'{e}or\`{e}me 2.4.2]{I}, \cite{G-S}).
Moreover, we have inclusions of \'{e}tale sheaves
\begin{equation*}
\lambda^{n}_{X, r}
\subset
W_{r}\omega^{n}_{X, \log}
\subset
\nu^{n}_{X, r}
\end{equation*}
(\cite[p.736, Proposition 4.2.1]{SaL}) and
these inclusions are not equalities in general 
(cf. \cite[p.737, Remark 4.2.3]{SaL}). 
In \S\ref{LH},
we prove the following:
\begin{Theorem}\upshape(Gersten resolution, Theorem \ref{ninj})
Let $A$ be a local ring of
a normal crossing variety over 
a field of positive characteristic $p>0$. 
Then the sequence
\begin{equation*}
0\to    
\HO^{s}_{\et}(A, \nu^{n}_{r})
\to
\bigoplus_{x\in \spec(A)^{(0)}}
\HO^{s}_{x}(A_{\et}, \nu^{n}_{r})
\to
\bigoplus_{x\in \spec(A)^{(1)}}
\HO^{s+1}_{x}(A_{\et}, \nu^{n}_{r})
\to
\cdots
\end{equation*}
is exact for any integers $s$ and $r>0$.
Here
\begin{math}
\HO^{t}_{x}(A_{\et}, \nu^{n}_{r})
=
\displaystyle
\lim
_{\substack{\to\\U}}
\HO^{t}_{\bar{\{x\}}\cap U}
(U_{\et}, \nu^{n}_{r})
\end{math}
for a non-negative integer $t$
where $U$ runs through open subschemes of
$\spec(A)$ such that
$x\in U$.
\end{Theorem}
\begin{Theorem}\upshape\label{Me1}%
(Rigidity, Theorem \ref{LL})
Let 
$X$ be a normal crossing variety over a field of characteristic $p>0$,
$A$ the henselization of a local ring 
$\mathcal{O}_{X, x}$ of $x$ in $X$
and $k$ the residue field of $A$.
Let $n$ be a non negative integer and $r$ a positive integer.
Then the homomorphism
\begin{equation}\label{ier}
\HO^{1}_{\et}\left(
A, \lambda^{n}_{A, r}
\right)
\to
\HO^{1}_{\et}\left(
k, \lambda^{n}_{k, r}
\right)
\end{equation}
is an isomorphism.
\end{Theorem}
In the case where $X$
is a smooth variety over a field
of characteristic $p>0$, 
Theorem \ref{Me1} is 
\cite[p.55, Theorem 5.3]{SM}.
The proof of
Theorem \ref{Me1} is reduced to 
\cite[p.55, Theorem 5.3]{SM}
by using the contravariant functoriality of the sheaves 
$\lambda^{n}_{X, r}$
for normal crossing varieties $X$ (\cite[p.734, Corollary 3.5.3]{SaL}).
If we replace $\lambda_{A, r}^{n}$
with
$\nu^{n}_{A, r}$ in the homomorphism (\ref{ier}),
the homomorphism (\ref{ier}) is not an isomorphism
unless $A$ is smooth over a field of positive characteristic 
(see Remark \ref{nrns} below).

\vspace{2.0mm}

Next, we consider mixed characteristic cases.
Let $B$ be a discrete valuation ring of
mixed characteristic $(0, p)$
with the quotient field $K$.
Let $\mathfrak{X}$ be a \textit{semistable family} over $\spec (B)$,
i.e. a regular scheme of pure dimension which is
flat of finite type over $\spec (B)$, 
$\mathfrak{X}\otimes_{B}K$ is smooth over $\spec (K)$,
and 
the special fiber $Y$ of $\mathfrak{X}$ is 
a reduced divisor with normal crossings on $\mathfrak{X}$.

Let $j$ and $\iota$ be as follows:
\begin{equation*}
\mathfrak{X}_{K}
\xrightarrow{j}  
\mathfrak{X}
\xleftarrow{\iota}
Y.
\end{equation*}
Then there is an exact sequence of sheaves on $\mathfrak{X}_{\et}$
\begin{equation*}
R^{n}j_{*}\mu_{p^{r}}^{\otimes n}
\to
\bigoplus_{y\in Y^{(0)}}i_{y*}W_{r}\Omega^{n-1}_{y, \log}
\to
\bigoplus_{y\in Y^{(1)}}i_{y*}W_{r}\Omega^{n-2}_{y, \log}
\end{equation*}
where 
$\mu_{p^{r}}$ is the sheaf of $p^{r}$-th roots of unity and
each arrow arises from the boundary maps of Galois cohomology
groups (cf.\cite[pp.522--523, Lemma 1.3.1.(1)]{SaP}). Hence we have the morphism
\begin{equation}\label{jti}
R^{n}j_{*}\mu_{p^{r}}^{\otimes n}
\to
\iota_{*}\nu_{Y, r}^{n-1}
\end{equation}
by the definition of $\nu_{Y, r}^{n-1}$. In this situation, the $p$-adic 
Tate twist $\mathfrak{T}_{r}(n)$ is defined by K.Sato as follows:
\begin{Definition}\upshape\label{DefpTt}
(\textit{$p$-adic Tate twist} $\mathfrak{T}_{r}(n)$, 
cf.\cite[p.537, Definition 4.2.4]{SaP})
Let the notations be the same as above. For $n=0$,
\begin{equation*}
\mathfrak{T}_{r}(0)_{\mathfrak{X}}
:=
\mathbb{Z}/p^{r}\mathbb{Z}.
\end{equation*}
For $n\geq 1$, $\mathfrak{T}_{r}(n)$ be defined as a complex which is fitted into
the following distinguished triangle
\begin{align*}
\iota_{*}\nu_{Y, r}^{n-1}[-n-1]
\to
\mathfrak{T}_{r}(n)_{\mathfrak{X}}
\to    
\tau_{\leq n}Rj_{*}\mu_{p^{r}}^{\otimes n}
\xrightarrow{\sigma_{\mathfrak{X}, r}(n)}
\iota_{*}\nu_{Y, r}^{n-1}[-n]
\end{align*}
where the morphism $\sigma_{\mathfrak{X}, r}(n)$ is induced 
by the morphism (\ref{jti}).
\end{Definition}

The main objective of this paper is to study 
\'{e}tale hypercohomology groups with values in $\mathfrak{T}_{r}(n)$.
By Theorem \ref{QH} and Theorem \ref{QS}
which are results about the $p$-adic vanishing cycle 
$\iota^{*}R^{n}j_{*}\mu_{p}^{\otimes n}$
due to Bloch-Kato/Hyodo and K.Sato, 
we have a finite filtration of 
$\mathcal{H}^{n}\left(\iota^{*}\mathfrak{T}_{1}(n)_{\mathfrak{X}}\right)$
which relates 
to the logarithmic Hodge-Witt sheaves $\lambda^{n}_{Y, 1}$
and the modified differential modules
(see Remark \ref{filtiT} below).
Then we are able to prove the following 
by using Theorem \ref{Me1} and an isomorphism (\ref{isohi}):
\begin{Theorem}\upshape\label{Me2}
(Rigidity, Theorem \ref{MR})
Let $\mathfrak{X}$ be the same as above,
$R$ the henselization of 
a local ring $\mathcal{O}_{\mathfrak{X}, x}$
of $x$ in $\mathfrak{X}$
and $k$ the residue field of $R$.
Let $n$ be a non negative integer and $r$
be
a positive integer.
Then we have an isomorphism
\begin{equation*}
\HO_{\et}^{n+1}(R, \mathfrak{T}_{r}(n))   
\simeq
\HO^{1}_{\et}(k, \lambda_{k, r}^{n}).
\end{equation*}
%
\end{Theorem}
Theorem \ref{Me2} has several applications.
Let us denote Bloch's cycle complex
for the Zariski topology
(\cite{B}, \cite{L}) by $\mathbb{Z}(n)$.
For a positive integer $m$,
$\mathbb{Z}(n)\otimes\mathbb{Z}/m$
denotes by $\mathbb{Z}/m(n)$.
Thus, $\mathbb{Z}(n)_{\et}$ (resp. $\mathbb{Z}/m(n)_{\et}$) denotes 
the \'{e}tale sheafification of
$\mathbb{Z}(n)$ (resp. $\mathbb{Z}/m(n)$).

Then, the first application of Theorem \ref{Me2} is as follows:
\begin{Proposition}\upshape\label{IVAM}(Proposition \ref{VAM})
Let $R$ be the same as 
in Theorem \ref{Me2}. 
Then we have
\begin{equation*}
\HO^{i}_{\Zar}(
R, 
\mathbb{Z}/p^{r}(n)
)
=0
\end{equation*}
for $i=n+1$.
\end{Proposition}

By Proposition \ref{IVAM} and \cite[p.209, Remark 7.2]{SaR}, 
we have an isomorphism
\begin{equation}\label{ISaC}
\tau_{\leq n+1}(\mathbb{Z}/p^{r}(n)_{\et})
\simeq
\mathfrak{T}_{r}(n)
\end{equation}
in $D^{b}(\mathfrak{X}_{\et}, \mathbb{Z}/p^{r}\mathbb{Z})$
which is the derived category of bounded complexes of \'{e}tale 
$\mathbb{Z}/p^{r}\mathbb{Z}$-sheaves on $\mathfrak{X}$.
In \cite[p.524, Conjecture 1.4.1 (1)]{SaP}, it is conjectured that the truncation in
the isomorphism (\ref{ISaC}) is unnecessary.
If $\mathfrak{X}$ is smooth over $\spec(B)$, then
this conjecture holds true by \cite[p.786, Corollary 4.4]{Ge}.

Then,
by using the isomorphism (\ref{ISaC}),
the second application of Theorem \ref{Me2} is expressed as follows: 
\begin{Theorem}\label{GCp}\upshape
(Gersten resolusion, Theorem \ref{LGe} and (\ref{ISaC}))
Let 
$\mathfrak{X}$ be a semistable family
over the spectrum of a discrete valuation ring 
of mixed characteristic $(0, p)$ and
$R=\tilde{\mathcal{O}_{\mathfrak{X}, x}}$ 
the henselian local ring
of $x$ in $\mathfrak{X}$. 
Suppose that $\operatorname{dim}(R)=2$.
Then the sequence
%
\begin{equation}\label{expaT}
0\to 
\HO^{n+1}_{\et}(
R,
\mathbb{Z}/m(n)
)
\to
\HO^{n+1}_{\et}
(
k(R),
\mathbb{Z}/m(n)
)
\to
\bigoplus_{\substack{
\mathfrak{p}\in\spec (R)^{(1)}
\\
}
}
\HO^{n+2}_{\mathfrak{p}}
(
(R_{\mathfrak{p}})_{\et},
\mathbb{Z}/m(n)
)
\end{equation}
%
is exact for any integer $n\geq 0$ and $m=p^{r}$. 
Here $k(R)$ is the fraction field of $R$.
\end{Theorem}
Let $B$ be a regular local ring of dimension at most $1$.
Let $l$ be a positive integer which is invertible in $B$
and $\mu_{l}$ the \'{e}tale sheaf of $l$-th roots of unity. 
Then we have an isomorphism
\begin{equation}\label{IBGl}
\tau_{\leq n+1}(\mathbb{Z}/l(n)_{\et})
\simeq
\mu_{l}^{\otimes n}
\end{equation}
in $D^{b}(\mathfrak{X}_{\et}, \mathbb{Z}/l\mathbb{Z})$ 
(see Remark \ref{BGl} below).
If 
$\mathfrak{X}$ is smooth over $\spec(B)$, 
then the truncation in (\ref{IBGl}) is unnecessary by 
\cite[Theorem 1.5]{G-L2} and \cite[p.786, Corollary 4.4]{Ge}. 
In the case where $R$ is equi-characteristic and 
$(m, \operatorname{char}(R))=1$, 
the exactness of the sequence (\ref{expaT})
is a part of the Bloch-Ogus Theorem (\cite{B-O})
by (\ref{IBGl}).

Moreover, we are able to prove the following local-global principle
as an application of Theorem \ref{GCp} (and Theorem \ref{Me2}):
\begin{Theorem}\label{Me3}\upshape
(Corollary \ref{LGLC} and Theorem \ref{MTP})
Let $B$ be an excellent henselian discrete valuation ring of mixed characteristic
$(0, p)$.
Let $X$ be a regular scheme over  $\spec(B)$  and $Y$ the special fiber of $X$.
Let $\tilde{\mathcal{O}_{X, \mathfrak{p}}}$ 
be the henselization of a local ring
$\mathcal{O}_{X, \mathfrak{p}}$ of $\mathfrak{p}$
in $X$,
$\kappa(\mathfrak{p})$ the residue field of
$\mathfrak{p}\in X$
and $k(X)$ (resp. 
$k(\tilde{\mathcal{O}_{X, \mathfrak{p}}})$) 
the ring of rational functions on $X$ (resp. 
$\spec(\tilde{\mathcal{O}_{X, \mathfrak{p}}})$).
Suppose that 
$\operatorname{dim}(X)=2$.

Then the natural map
\begin{equation*}
\HO^{n+1}_{\et}
(
k(X),
\mu_{p^{r}}^{\otimes n}
)
\to
\bigoplus_{\mathfrak{p}\in 
X^{(1)}\setminus Y^{(0)}
}
\HO^{n}_{\et}
(
\kappa(\mathfrak{p}),
\mu_{p^{r}}^{\otimes (n-1)}
)
\oplus
\bigoplus_{\mathfrak{p}\in Y^{(0)}
}
\HO^{n+1}_{\et}
(
k(\tilde{\mathcal{O}_{X, \mathfrak{p}}})
,
\mu_{p^{r}}^{\otimes n}
)
\end{equation*}
%
is injective for $n\geq 1$
in the following cases:
\begin{description}
\item[(i)] $X$ is a scheme $\mathfrak{X}$ which is
a proper and semistable family 
over $\spec(B)$.
\item[(ii)] $X=\tilde{\mathcal{O}_{\mathfrak{X}, \mathfrak{p}}}$.
Here $\mathfrak{X}$ is the same as in (i) and 
$\mathfrak{p}\in \mathfrak{X}^{(2)}$.
\end{description}

This implies that the local-global map
\begin{equation}\label{GLMI}
\HO^{n+1}_{\et}
(
k(X),
\mu_{m}^{\otimes n}
)
\to
\prod_{\mathfrak{p}\in X^{(1)}}
\HO^{n+1}_{\et}
(
k(\tilde{\mathcal{O}_{X, \mathfrak{p}}})
,
\mu_{m}^{\otimes n}
)
\end{equation}
%
is injective for $n\geq 1$ and $m=p^{r}$.
\end{Theorem}
Suppose that $B$ is equi-characteristic and 
$(m, \operatorname{char}(B))=1$.
In the case (i),
Colliot-Th\'{e}l\`{e}ne raised a question whether 
the local-global map (\ref{GLMI}) is injective (\cite{C}).
In \cite[p.245, Theorem 3.3.6]{H-H-K}, 
Harbater–Hartmann-Krashen provided an affirmative answer
to this question
in this case.

In the case where $B$ is mixed characteristic $(0, p)$,
due to Hu's method 
(cf.\cite[The proof of Theorem 2.5]{Hu}, \cite[Remark 2.6 (2)]{Hu}),
Theorem \ref{Me3} (i) follows from 
Theorem \ref{GCp}
(or Theorem \ref{Me3} (ii))
without assuming $(m, p)=1$.
If $(m, p)=1$, 
then 
the exactness of the sequence (\ref{expaT}) has been proved in
\cite[p.34, Theorem 1]{SM19}. Theorem \ref{Me3} (i) is an extension of 
\cite[pp.62--63, Theorem 6.3]{SM}.
%

\section*{Acknowledgments}
I thank
Indian Institute of Science Education and Research, Bhopal
where the main part of this work was completed.

\subsection*{Notations}
For a scheme $X$,
$X^{(i)}$ denotes the set of points on $X$
of codimension $i$. 
Moreover, $k(X)$ 
denotes the ring of
rational functions on $X$ and $\kappa(x)$ 
denotes the residue field of $x\in X$.
If $X=\spec(A)$ for a ring $A$, 
$k(A)$ denotes $k(\spec (A))$.
For a scheme $X$, $X_{\Zar}$, $X_{\Nis}$ 
and $X_{\et}$
denote the category of \'{e}tale schemes over $X$
equipped with the Zariski, Nisnevich and \'{e}tale topology,
respectively.
For a scheme $X$, 
$D(X)$ denotes the derived category of complexes
of \'{e}tale sheaves of abelian groups on $X$, and
$D_{+}(X)\subset D(X)$ denotes the full subcategory
of complexes that are bounded below.
\section{
Logarithmic Hodge-Witt sheaves
}\label{LH}
In this section, we study the logarithmic Hodge-Witt sheaves
$\nu_{X, r}^{n}$ and $\lambda_{X, r}^{n}$ (\cite{SaL})
for a normal crossing variety $X$.
\subsection{$\nu^{n}_{X, r}$}
Let $X$ be a normal crossing variety over a field
of characteristic $p>0$.
Let $n$ be a non negative integer 
and $r$ a positive integer. 
We show a relation between two Zariski sheaves on $X_{\Zar}$
\begin{equation}\label{defzn}
\nu^{n}_{X, r}:=\operatorname{Ker}
\left(
\partial: 
\bigoplus_{x\in X^{(0)}}
i_{x *}W_{r}\Omega^{n}_{x, \log}
\to
\bigoplus_{x\in X^{(1)}}
i_{x *}W_{r}\Omega^{n-1}_{x, \log}
\right)
\end{equation}
(which is the image of the \'{e}tale sheaf
$\nu^{n}_{X, r}$ in \cite[p.527, 2.2]{SaP}
under the forgetful functor) 
and $\mathbb{Z}/p^{r}(n)$, where
$\mathbb{Z}(n)$ is
the Bloch’s cycle complex for the Zariski
topology (\cite{B},\cite{L}) and 
\begin{math}
\mathbb{Z}/p^{r}(n)
=\mathbb{Z}(n)\otimes\mathbb{Z}/p^{r}.
\end{math}
\begin{Proposition}\upshape\label{van}
Let $X$ be a normal crossing variety over a field.
Then we have
\begin{equation*}
\mathcal{H}^{i}
\left(
\mathbb{Z}(n)
\right)
=0
\end{equation*}
for $i\geq n+1$.
\end{Proposition}
\begin{proof}\upshape
It suffices to show the statement in the case where
$X$ is the spectrum of a local ring of a normal 
crossing variety. 
Moreover, the stalk of the sheaf $\mathcal{H}^{i}(\mathbb{Z}(n))$
at a point $x\in X$ is equal to
\begin{math}
\HO^{i}_{\Zar}(
\mathcal{O}_{X, x},
(j_{x})^{*}\mathbb{Z}(n)
)
\end{math}
by the spectral sequence
\begin{equation*}
E^{s, t}_{2}=\HO^{s}_{\Zar}(
\mathcal{O}_{X, x},
\mathcal{H}^{t}(\mathbb{Z}(n))
)
\Rightarrow
E^{s+t}=
\HO^{s+t}_{\Zar}(
\mathcal{O}_{X, x},
\mathbb{Z}(n)
)
\end{equation*}
and the definition of the Zariski cohomology.
Here $j_{x}: \spec(\mathcal{O}_{X, x})\to X$ 
is the natural map.
Hence it suffices to show
\begin{equation}\label{LV}
\HO^{i}_{\Zar}(\spec(A), \mathbb{Z}(n))  
=0
\end{equation}
for $i\geq n+1$
in the case where $A$ is a local ring of a normal crossing variety. 
We prove the equation (\ref{LV}) by induction on 
$\#\left(\spec (A)^{(0)}\right)$.

In the case where
$\#\left(\spec (A)^{(0)}\right)=1$,
the equation (\ref{LV}) follows from
\cite[p.786, Corollary 4.4]{Ge}.

Assume that the equation (\ref{LV}) holds
in the case where
$\#\left(\spec (A)^{(0)}\right)\leq r$.
Let 
\begin{equation}\label{Sr}
\spec (A_{1}), \cdots, \spec (A_{r+1})
\end{equation}
be the irreducible components of
\begin{math}
\displaystyle\spec (A).
\end{math}
Then elements of the set (\ref{Sr})
and
\begin{equation*}
\spec (B_{1})
=
\spec (A_{1})\bigcap\spec (A_{r+1}),
\cdots,
\spec (B_{r})
=
\spec (A_{r})\bigcap\spec (A_{r+1})
\end{equation*}
are smooth varieties by the definition of normal crossing variety.
Let
\begin{equation*}
\spec (B)
=
\displaystyle
\bigcup_{j=1}^{j=r}
\spec (B_{j})
=
\left(\displaystyle
\bigcup_{j=1}^{j=r}
\spec (A_{j})
\right)\cap
\spec (A_{r+1})
.
\end{equation*}
Then we have
\begin{equation}\label{BV}
\HO^{i}_{\Zar}\left(
\spec(B), \mathbb{Z}(n)
\right)
=0
\end{equation}
for $i\geq n+1$
by the inductive hypothesis.
Since 
we have
\begin{equation*}
\HO^{n+1}_{\Zar}
(\spec (A_{r+1})
,
\mathbb{Z}(n)
)
=0  
\end{equation*}
by \cite[p.786, Corollary 4.4]{Ge},
the homomorphism
\begin{equation*}
\HO_{\Zar}^{n}\left(
\spec (A_{r+1})\setminus\spec (B),
\mathbb{Z}(n)
\right)
\to
\HO_{\Zar}^{n-1}
\left(
\spec (B), \mathbb{Z}(n-1)
\right)
\end{equation*}
is surjective by the localization theorem 
(
\cite[p.277, Theorem (3.1)]{B}, \cite[p.537, Corollary (0.2)]{B2})
).
Moreover we have
\begin{equation*}
\spec (A)\setminus\spec (B)
=
(
\displaystyle
\bigcup_{j=1}^{j=r}
\spec (A_{j})
\setminus\spec (B)
)\bigoplus
\left(
\spec (A_{r+1})\setminus\spec (B)
\right)
\end{equation*}
and a commutative diagram
\begin{equation*}
\xymatrix{
\HO^{n}_{\Zar}(\spec(A_{r+1})\setminus\spec(B), 
\mathbb{Z}(n))
\ar[r]\ar[d]
&
\ar@{=}[d]  
\HO^{n-1}_{\Zar}(
\spec(B), \mathbb{Z}(n-1))
\\
\HO^{n}_{\Zar}(\spec(A)\setminus\spec(B), 
\mathbb{Z}(n))
\ar[r]
&
\HO^{n-1}_{\Zar}(
\spec(B), \mathbb{Z}(n-1)).
}
\end{equation*}
So the homomorphism 
\begin{equation}\label{ABS}
\HO_{\Zar}^{n}\left(
\spec (A)\setminus\spec (B),
\mathbb{Z}(n)
\right)
\to
\HO_{\Zar}^{n-1}
\left(
\spec (B), \mathbb{Z}(n-1)
\right)
\end{equation}
is also surjective. Since
we have the equation (\ref{BV}) and
$A_{j}$ $(1\leq j\leq r+1)$ are local rings of a 
smooth variety,
we have
\begin{equation*}
\HO_{\Zar}^{i}
(
\displaystyle
\bigcup_{j=1}^{j=r}
\spec (A_{j})
\setminus\spec (B)
, 
\mathbb{Z}(n)
)
=
\HO_{\Zar}^{i}
\left(
\spec (A_{r+1})
\setminus\spec (B)
, \mathbb{Z}(n)
\right)
=0
\end{equation*}
for $i\geq n+1$ by the inductive hypothesis and the localization theorem.
Hence we have
\begin{equation}\label{ABV}
\HO_{\Zar}^{i}
\left(
\spec (A)\setminus\spec (B), 
\mathbb{Z}(n)
\right)
=0   
\end{equation}
for $i\geq n+1$. By using the localization theorem, the sequence
%
\begin{align*}
&\HO^{i}_{\Zar}(\spec(A)\setminus\spec(B), \mathbb{Z}(n))
\to
\HO^{i-1}_{\Zar}(\spec(B), \mathbb{Z}(n-1)) \\
\to
&\HO^{i+1}_{\Zar}(\spec(A), \mathbb{Z}(n)) 
\to
\HO^{i+1}_{\Zar}(\spec(A)\setminus\spec(B), \mathbb{Z}(n))
\end{align*}
%
is exact for any integer $i$.
Hence the equation (\ref{LV}) follows from 
the equations (\ref{BV}), 
(\ref{ABV}) and 
the surjectivity of the homomorphism (\ref{ABS}). 
This completes the proof.
\end{proof}
\begin{Proposition}\upshape\label{apbk}
Let $X$ be a normal crossing variety over a field.
\begin{itemize}
\item[1.]
If $\iota: Z\to X$ is a closed subscheme of codimension $c$, 
then the canonical map
\begin{equation*}
\tau_{\leq n+2}
\left(
\mathbb{Z}(n-c)_{\et}^{Z}[-2c]
\right)
\to
\tau_{\leq n+2}R\iota^{!}
\mathbb{Z}(n)_{\et}^{X}
\end{equation*}
is a quasi-isomorphism, where 
$\mathbb{Z}(n)_{\et}$
is the \'{e}tale sheafification of
$\mathbb{Z}(n)$.
\item[2.]
Let 
$\epsilon: X_{\et}\to X_{\Zar}$
be the change of site map.
Then the canonical map induces a quasi-isomomorphism
\begin{equation*}
\mathbb{Z}(n)_{\Zar}
\xrightarrow{\sim}
\tau_{\leq n+1}
R\epsilon_{*}\mathbb{Z}(n)_{\et}.
\end{equation*}
\end{itemize}
\end{Proposition}
\begin{proof}\upshape
By the same argument as in the proof of \cite[Proposition 2.1]{SM}, 
Proposition \ref{apbk}.1
follows from 
Proposition \ref{apbk}.2.
So it suffices to show Proposition \ref{apbk}.2.
Moreover, it suffices to show Proposition \ref{apbk}.2 in the case where $X$
is the spectrum of a local ring $A$ of 
a normal crossing variety.

We prove Proposition $\ref{apbk}.2$ by induction on $\#(\spec(A)^{(0)})$. 
Suppose that 
$\#(\spec(A)^{(0)})=1$.
Then we have a quasi-isomorphism
\begin{equation*}
\mathbb{Q}/\mathbb{Z}(n)_{\Zar}
\xrightarrow{\sim}
\tau_{\leq n}R\epsilon_{*}
\mathbb{Q}/\mathbb{Z}(n)_{\et}
\end{equation*}
by \cite[Theorem 8.5]{G-L}, \cite[Theorem 1.6]{G-L2} and \cite{V2}. 
Moreover we have a quasi-isomorphism
\begin{equation*}
\mathbb{Q}(n)_{\Zar}  
\xrightarrow{\sim}
R\epsilon_{*}\mathbb{Q}(n)_{\et}
\end{equation*}
by \cite[p.781, Proposition 3.6]{Ge}.
So Proposition \ref{apbk}.2 holds in the case where
$\#(\spec(A)^{(0)})=1$ by the five lemma. 
See also \cite[p.774, Theorem 1.2.2]{Ge}.

Assume that Proposition \ref{apbk}.2 holds in the
case where $\#(\spec(A)^{(0)})\leq r$.
Suppose that $\#(\spec(A)^{(0)})=r+1$
and $\spec(A_{k})~(1\leq k\leq r+1)$ are
the irreducible components of
\begin{math}
\displaystyle\spec (A).
\end{math}

Let
$\spec(B)=\spec(A)\cap\spec(A_{r+1})$ and $\iota: \spec(B)\to \spec(A)$ be 
the closed immersion. Let 
$j: U\to \spec(A)$ be the open complement of $\iota$. Then we have
\begin{equation*}
U=  
(
\displaystyle
\bigcup_{k=1}^{k=r}
\spec (A_{k})
\setminus\spec (B)
)\bigoplus
\left(
\spec (A_{r+1})\setminus\spec (B)
\right).
\end{equation*}
Moreover,
\begin{math}
\displaystyle
\bigcup_{k=1}^{k=r}
\spec (A_{k})
\setminus\spec (B)
\end{math}
and
\begin{math}
\spec(A_{r+1})
\setminus
\spec(B)
\end{math}
are normal crossing varieties,
\begin{equation*}
\#\left(
\left(
\displaystyle
\bigcup_{k=1}^{k=r}
\spec (A_{k})
\setminus\spec (B)
\right)^{(0)}
\right)
=r
~~
\textrm{and}
~~
\#\left(
\left(
\spec(A_{r+1})
\setminus
\spec(B)
\right)^{(0)}
\right)
=1.
\end{equation*}
Hence Proposition \ref{apbk}.2 holds for $U$ by the assumption.
Then Proposition \ref{apbk}.1 holds for 
$\iota: \spec(B)\to\spec(A)$
by Proposition \ref{van} and
the same argument as in the proof of 
\cite[Proposition 2.1]{SM}.
Moreover, 
$\#(\spec(B)^{(0)})=r$ and
Proposition \ref{apbk}.2 holds for
$\spec(B)$ by the assumption.
Therefore Proposition \ref{apbk}.2 follows by
Proposition \ref{van} and
the same argument as in the proof of \cite[Theorem 1.2.2]{Ge}.
This completes the proof.
\end{proof}
\begin{Proposition}\upshape\label{iso}
Let $X$ be a normal crossing variety over a field
of characteristic $p>0$. 
Then
we have an isomorphism in 
$D^{b}(X_{\Zar}, \mathbb{Z}/p^{r}\mathbb{Z})$
\begin{equation*}
\mathbb{Z}/p^{r}(n)
\simeq
\nu^{n}_{X, r}[-n],
\end{equation*}
where 
$D^{b}(X_{\Zar}, \mathbb{Z}/p^{r}\mathbb{Z})$
is
the derived category of 
bounded complexes of
Zariski 
$\mathbb{Z}/p^{r}\mathbb{Z}$ -sheaves on $X$.
\end{Proposition}
\begin{proof}\upshape
It suffices to show the statement in the case where
$X$ is the spectrum of a local ring of a normal crossing variety.
We prove the statement by induction on 
$\operatorname{dim}X$. In the case where 
$\operatorname{dim}X=0$, the statement follows from 
\cite[p.491, Theorem 8.3]{G-L}.

Assume that the statement holds in the case where
$\operatorname{dim}X\leq m$. Suppose that $\operatorname{dim}X=m+1$.
Let $\iota: Y\to X$ be a closed immersion of codimension $1$,
$Y$ regular and 
$j: U\to X$ the complement of $\iota$.
Let $Y$ be the spectrum of a regular local ring $A^{\prime}$. 
Then $A^{\prime}$ is a local ring of a regular
ring of finite type over a field of positive characteristic. 
By Quillen's method 
(cf.\cite[\S7, The proof of Theorem 5.11]{Q}),
\begin{equation}\label{QM}
A^{\prime}
=
\displaystyle
\lim_{\to}
A_{s}^{\prime}
\end{equation}
where $A_{s}^{\prime}$ is a local ring of a smooth algebra over 
$\mathbb{F}_{p}$ and the maps $A_{s}^{\prime}\to A_{t}^{\prime}$
are flat.
So Proposition \ref{iso} holds for $Y$ by \cite[p.491, Theorem 8.3]{G-L}.
Since
$\operatorname{dim}U=m$
and the sequence
\begin{equation*}
\cdots
\to
\iota_{*}\mathbb{Z}/p^{r}(n-1)^{Y}[-2]
\to
\mathbb{Z}/p^{r}(n)^{X}
\to
Rj_{*}\mathbb{Z}/p^{r}(n)^{U}
\to
\cdots
\end{equation*}
is a distinguished triangle by \cite[p.277, Theorem (3.1)]{B}, 
we have a quasi-isomorphism
\begin{equation*}
\tau_{\leq n}
\left(
\mathbb{Z}/p^{r}(n)
\right)
\simeq
\nu_{X, r}^{n}[-n]
\end{equation*}
by (\ref{defzn}).
Moreover, we have
\begin{equation*}
\mathcal{H}^{s}
\left(
\mathbb{Z}/p^{r}(n)
\right)
=0
\end{equation*}
for $s>n$ by Proposition \ref{van}. This completes the proof.
\end{proof}
\begin{Corollary}\upshape\label{GeZ}
Let $A$ be a local ring at a point 
of a normal crossing variety over 
a field of positive characteristic $p>0$. 
Then the sequence
\begin{align*}
0\to \HO^{0}(A, \nu^{n}_{r})
\to
\bigoplus_{x\in \spec(A)^{(0)}}
\HO^{0}(\kappa(x), \nu_{r}^{n})
\to
\bigoplus_{x\in \spec(A)^{(1)}}
\HO^{0}(\kappa(x), \nu_{r}^{n-1})  \\
\to
\bigoplus_{x\in \spec(A)^{(2)}}
\HO^{0}(\kappa(x), \nu_{r}^{n-2})
\to
\cdots
\to
\bigoplus_{x\in \spec(A)^{(i)}}
\HO^{0}(\kappa(x), \nu_{r}^{n-i})
\to\cdots
\end{align*}
is exact for any integer $n$.
\end{Corollary}
\begin{proof}\upshape
We consider the spectral sequence
\begin{equation}\label{sph}
E^{s, t}_{1}=
\bigoplus_{x\in \spec(A)^{(s)}}
\HO^{2n-s+t}_{\Zar}(
\kappa(x), \mathbb{Z}/p^{r}(n-s))
\Rightarrow
E^{s+t}=
\HO^{2n+s+t}_{\Zar}(A, \mathbb{Z}/p^{r}(n))    
\end{equation}
(cf.\cite[p.782]{Ge}). 
In \cite{Ge} the spectral sequence (\ref{sph}) is shown in the smooth case.
But the construction of the spectral sequence
(\ref{sph}) can be carried out without assuming smoothness. 
By Proposition \ref{iso}, we have
$E^{s, t}_{1}=0$ for $t\neq -n$ and
$E^{s+t}=0$ for $s+t\neq -n$.
Hence the statement follows from the spectral sequence (\ref{sph}).
\end{proof}
\begin{Remark}\upshape
If $A$ is the strict henselization of a local ring of
a normal crossing variety over a field of positive characteristic,
Corollary \ref{GeZ} directly follows from 
\cite[p.716, Corollary 2.2.5 (1)]{SaL}.
\end{Remark}
\begin{Proposition}\upshape\label{Sm2v}
Let $A$ be an equidimensional catenary local ring of characteristic
$p>0$
and 
$F$ be an  \'{e}tale 
$p$-torsion sheaf.

Then we have
\begin{equation*}
\HO^{r+j}_{x}(A_{\et}, F)=0    
\end{equation*}
for $r\geq 2$ and $x\in \spec(A)^{(j)}$.
\end{Proposition}
\begin{proof}\upshape
Let $A_{x}$ be the local ring of $A$ at
$x\in\spec(A)$.
Then
\begin{math}
\operatorname{codim}\left(\bar{\{x\}}, \spec(A)\right)=\operatorname{dim}(A_{x})
\end{math}
and we have
\begin{equation*}
\HO^{s}_{x}(A_{\et}, F)
=\displaystyle\lim_{\substack{
\to  \\
U}
}
\HO^{s}_{\bar{\{x\}}\cap U}(U_{\et}, F)
=\HO^{s}_{x}((A_{x})_{\et}, F)
\end{equation*}
for any integer $s\geq 0$ 
by \cite[pp.88--89, III, Lemma 1.16]{M} 
where $U$ runs through open subschemes of $\spec(A)$ such that $x\in U$.
So it suffices to show the statement in the case where
$x\in \spec(A)$ is the closed point of $\spec(A)$.
We prove the statement by induction on 
\begin{math}
\operatorname{dim}(A).
\end{math}

In the case where $\operatorname{dim}(A)=0$,
the statement is true because
$A$ is a ring of characteristic $p>0$
and $p$-cohomological dimension 
$\operatorname{cd}_{p}(A)\leq 1$.

Assume that the statement is true for 
$\operatorname{dim}(A_{x})\leq j$. 
Then we prove the statement in the case where 
$\operatorname{dim}(A)=j+1$.
In order to prove the statement,
we use the spectral sequence
\begin{equation}\label{spp}
E^{s, t}_{1}
=
\displaystyle
\bigoplus_{x\in \spec(A)^{(s)}}
\HO^{s+t}_{x}(A_{\et}, F)
\Rightarrow
E^{s+t}
=
\HO^{s+t}_{\et}(A, F)
\end{equation}
(cf.\cite[Part 1, \S 1]{C-H-K}).
Since 
\begin{math}
\operatorname{dim}(A)=j+1,
\end{math}
we have
\begin{equation*}
E^{s, t}_{1}=0    
\end{equation*}
for $s>j+1$. So we have
\begin{equation}\label{VE2}
E^{s, t}_{2}=0    
\end{equation}
for $s>j+1$.

Moreover, we have (\ref{VE2}) for $s\leq j$
and $t\geq 2$ by the assumption.
Hence we have
\begin{equation}\label{E2inf}
E_{2}^{j+1, r}
=
E_{\infty}^{j+1, r}
\end{equation}
for $r\geq 2$. 
Since 
$\operatorname{cd}_{p}(A)\leq 1$ by 
\cite[Expos\'{e} X, Th\'{e}or\`{e}me 5.1]{SGA4}, 
we have
\begin{equation}\label{infv}
E_{\infty}^{j+1, r}=E^{j+r+1}
=0    
\end{equation}
for $r\geq 2$ by the spectral sequence (\ref{spp}). 
Hence the sequence
\begin{equation*}
\displaystyle
\bigoplus_{x\in \spec(A)^{(j)}}
\HO_{x}^{r+j}(A_{\et}, F)
\to
\displaystyle
\bigoplus_{x\in \spec(A)^{(j+1)}}
\HO_{x}^{r+j+1}(A_{\et}, F)
\to 0    
\end{equation*}
is exact for $r\geq 2$
by (\ref{E2inf}) and (\ref{infv}).
Therefore we have
\begin{equation*}
\HO_{x}^{r+j+1}(A_{\et}, F)=0    
\end{equation*}
for $r\geq 2$ by the assumption. This completes the proof.
\end{proof}
\begin{Theorem}\upshape\label{ninj}
Let $A$ be a local ring of
a normal crossing variety over 
a field of positive characteristic $p>0$. 
Then the sequence
\begin{equation*}
0\to    
\HO^{t}_{\et}(A, \nu^{n}_{r})
\to
\bigoplus_{x\in \spec(A)^{(0)}}
\HO^{t}_{x}(A_{\et}, \nu^{n}_{r})
\to
\bigoplus_{x\in \spec(A)^{(1)}}
\HO^{t+1}_{x}(A_{\et}, \nu^{n}_{r})
\to
\cdots
\end{equation*}
is exact for any integers $t$ and $r>0$.
\end{Theorem}
\begin{proof}\upshape
In the case where $t<0$, the statement follows from \cite[p.718, Theorem 2.4.2]{SaL}.
In the case where $t=0$, the statement follows from 
Corollary \ref{GeZ} and \cite[p.718, Theorem 2.4.2]{SaL}.
In the case where $t\geq 2$, the statement follows from
Proposition \ref{Sm2v} and \cite[Expos\'{e} X, Th\'{e}or\`{e}me 5.1]{SGA4}.
So it suffices to show the statement in the case where $t=1$.
In order to prove the statement,
we use the spectral sequence (\ref{spp}) for $F=\nu_{r}^{n}$.
By Corollary \ref{GeZ} and \cite[p.718, Theorem 2.4.2]{SaL},
we have
\begin{equation}\label{E2v2}
E_{2}^{s, t}=0    
\end{equation}
for $s>0$ and $t\leq 0$. By Proposition \ref{Sm2v}, 
we have the equation (\ref{E2v2}) for $t\geq 2$.
Moreover, we have 
\begin{math}
E^{s+t}=0
\end{math}
for $s+t\geq 2$ by \cite[Expos\'{e} X, Th\'{e}or\`{e}me 5.1]{SGA4}.
Hence we have
isomorphisms
\begin{align*}
E^{s, 1}_{2}
\simeq
\left\{
\begin{array}{ll}
E^{1} & (s=0)
\\
0  & (s>0)
\end{array}
\right.
\end{align*} 
by the spectral sequence (\ref{spp}) for $F=\nu^{n}_{r}$.
This completes the proof.
\end{proof}
\subsection{$\lambda_{X, r}^{n}$}
Let $n$ be a non negative integer and $r$ a positive integer.
For a normal crossing variety $X$ over a field of positive characteristic, 
we define 
a kind of generalized logarithmic Hodge-Witt sheaves on $X_{\et}$ by
\begin{equation*}
\lambda^{n}_{X, r}
:=
\operatorname{Im}\left(
d\log:
(\mathbb{G}_{m, X})^{\otimes n}
\to
\bigoplus_{x\in X^{(0)}}i_{x*}\operatorname{W}_{r}\Omega_{x, \log}^{n}
\right)
\end{equation*}
(cf.\cite[p.726, Definition 3.1.1]{SaL}). Then we show the following:

\begin{Theorem}\upshape\label{LL}
Let $A$ be the henselization of a local ring of 
a normal crossing variety over a field of characteristic $p>0$
and $k$ the residue field of $A$.
Then the homomorphism
\begin{equation}\label{Ru}
\HO^{1}_{\et}\left(
A, \lambda^{n}_{A, r}
\right)
\to
\HO^{1}_{\et}\left(
k, \lambda^{n}_{k, r}
\right)
\end{equation}
is an isomorphism.
\end{Theorem}
\begin{proof}\upshape
Let $C$ be a finitely generated $k$-algebra. 
Then there exists a surjective $k$-algebra homomorphism
\begin{equation*}
k[T_{1}, \cdots, T_{N}] 
\to
C
\end{equation*}
for some integer $N$.
Hence $A$ is embedded into an equi-characteristic henselian regular local ring $B$.
Let $\iota: \spec (A)\to \spec (B)$ be the closed immersion. Then
the natural pull-back map
\begin{equation}\label{Pul}
\iota^{*}: 
\mathbb{G}_{m, B}^{\otimes n}
\to
\iota_{*}\left(
\mathbb{G}_{m, A}^{\otimes n}
\right)
\end{equation}
is surjective. 
Moreover the homomorphism (\ref{Pul})
induces
the pull-back map
\begin{equation}\label{Pul2}
\iota^{*}:
\lambda_{B, r}^{n}
\to
\iota_{*}\lambda_{A, r}^{n}
\end{equation}
(cf.\cite[p.732, Theorem 3.5.1]{SaL})
and the homomorphism (\ref{Pul2})
is also surjective. 
Since $B$ has $p$-cohomological dimension at most $1$,
we have
\begin{equation*}
\HO^{2}_{\et}(
B,
\operatorname{Ker}(\iota^{*})
)
=0.
\end{equation*}
%
Hence
the homomorphism
\begin{equation}\label{BA}
\HO^{1}_{\et}(B, \lambda_{B, r}^{n})
\to
\HO^{1}_{\et}(A, \lambda_{A, r}^{n})
\end{equation}
which is induced by the pull-back map (\ref{Pul2})
is surjective. 
Moreover a composition 
of the homomorphism (\ref{BA}) and
the homomorphism
\begin{equation}\label{Ak}
\HO^{1}_{\et}(A, \lambda_{A, r}^{n})
\to
\HO^{1}_{\et}(k, \lambda_{k, r}^{n})   
\end{equation}
which is induced by the pull-back map
coincides with the homomorphism
\begin{equation}\label{Bk}
\HO^{1}_{\et}(B, \lambda_{B, r}^{n})
\to
\HO^{1}_{\et}(k, \lambda_{k, r}^{n})     
\end{equation}
which is induced by the pull-back map by 
\cite[p.734, Corollary 3.5.3]{SaL}.
Then the homomorphism
(\ref{Bk})
is an isomorphism by \cite[Theorem 5.3]{SM}. Therefore
the homomorphism (\ref{BA}) is an isomorphism and
the homomorphism (\ref{Ak}) is also an isomorphism. This completes
the proof.
\end{proof}
\begin{Remark}\upshape\label{nrns}
If we replace $\lambda_{A, r}^{n}$ 
with $\nu_{A, r}^{n}$ in the homomorphism (\ref{Ru}),
the homomorphism (\ref{Ru}) is not an isomorphism in general.
In fact,
\begin{align*}
\lambda_{A, r}^{0}=\mathbb{Z}/p^{r}
&&
\nu^{0}_{A, r}
=
\underbrace{
\mathbb{Z}/p^{r}\oplus\cdots\oplus\mathbb{Z}/p^{r}
}_{\#\left(\spec (A)^{(0)}\right)~
\textit{times}}
\end{align*}
and
\begin{equation*}
\HO^{1}_{\et}(A, \lambda_{A, r}^{0}) 
\neq
\HO^{1}_{\et}(A, \nu_{A, r}^{0}) 
\end{equation*}
unless $A$ is smooth over a field 
of positive characteristic.
\end{Remark}
\section{$p$-adic Tate twist $\mathfrak{T}_{r}(n)$}
Let $B$ be a discrete valuation ring of
mixed characteristic $(0, p)$ with the quotient field $K$.

Let
$\mathfrak{X}$ be a 
semistable family
over $\spec (B)$, 
i.e. 
a regular scheme of pure dimension
which is flat of finite type over $\spec (B)$, 
$\mathfrak{X}_{K}=\mathfrak{X}\otimes K$ is smooth over $\spec (K)$,
and the special fiber $Y$ of $\mathfrak{X}$
is a reduced divisor with normal crossings on $\mathfrak{X}$.

Let $j$ and $\iota$ be as follows:
\begin{equation*}
\mathfrak{X}_{K}\xrightarrow{j}
\mathfrak{X}  
\xleftarrow{\iota}
Y.
\end{equation*}
Let $n$ be a non negative integer
and $r$  a positive integer.
In this section, we study 
the $p$-adic Tate twist
$\mathfrak{T}_{r}(n)$ (cf.\cite[p.537, Definition 4.2.4]{SaP}).
In order to study it, 
the $p$-adic vanishing cycle
\begin{equation*}
M^{n}_{r}
:=
\iota^{*}R^{n}j_{*}\mu_{p^{r}}^{\otimes n}
\end{equation*}
plays important roles where $\mu_{p^{r}}$ is the sheaf of $p^{r}$-th
roots of unity.
\begin{Remark}\upshape
In \cite{SaP}, 
$\mathfrak{T}_{r}(n)$ is defined in the case where the residue field of $B$ is perfect.
But the results in \cite{SaP} hold true without this assumption as explained in
\cite[p.187, Remark 3.7]{SaR}.
\end{Remark}
\subsection{Review on the structure of the $p$-adic vanishing cycle}
Let the notations be the same as above.
We define the \'{e}tale sheaf 
$\mathcal{K}^{M}_{n, \mathfrak{X}_{K}/Y}$ on $Y$
as
\begin{equation*}
\mathcal{K}^{M}_{n, \mathfrak{X}_{K}/Y}
:=
(
\iota^{*}j_{*}
\mathcal{O}^{\times}_{\mathfrak{X}_{K}}
)^{\otimes n}
/J,
\end{equation*}
where $J$ denotes the subsheaf which is generated by
local sections of the form
\begin{equation*}
x_{1}\otimes \cdots\otimes x_{n} 
~(x_{i}\in \iota^{*}j_{*}\mathcal{O}^{\times}_{\mathfrak{X}_{K}})
~\textrm{with}~
x_{i}+x_{j}=0
~\textrm{or}~1
\end{equation*}
for some $1\leq i<j\leq n$.
By \cite[(1.2)]{B-K}, there is a natural map
\begin{equation}\label{MM}
\mathcal{K}^{M}_{n, \mathfrak{X}_{K}/Y}
\to
M^{n}_{r}
\end{equation}
and we define the filtrations $U^{*}$ and $V^{*}$ 
on 
the $p$-adic vanishing cycle $M^{n}_{r}$ 
by using this map (\ref{MM}), as follows:
\begin{Definition}\upshape(cf.\cite[p.546, (1.4)]{H}, \cite[pp.530--531, Definition 3.3.2]{SaP})
\begin{itemize}
\item[(1)] Let $\pi$ be a prime element of $B$. Let
$U^{0}_{\mathfrak{X}_{K}}$ be the full sheaf 
$\iota^{*}j_{*}\mathcal{O}^{\times}_{\mathfrak{X}_{K}}$.
For $q\geq 1$, let $U^{q}_{\mathfrak{X}_{K}}$ be the 
\'{e}tale subsheaf of 
$\iota^{*}j_{*}\mathcal{O}^{\times}_{\mathfrak{X}_{K}}$
which is generated by local sections of the
form $1+\pi^{q}\cdot a$ with 
$a\in\iota^{*}\mathcal{O}_{\mathfrak{X}}$.
We define the subsheaf 
$U^{q}\mathcal{K}^{M}_{n, \mathfrak{X}_{K}/Y}$
as the part which is generated by
\begin{math}
U^{q}_{\mathfrak{X}_{K}}
\otimes
\{
\iota^{*}j_{*}\mathcal{O}^{\times}_{\mathfrak{X}_{K}}
\}^{\otimes (n-1)}.
\end{math}
\item[(2)]
We define the subsheaf $U^{q}M^{n}_{r}$ ($q\geq 0$)
of $M^{n}_{r}$ as the image of
$U^{q}\mathcal{K}^{M}_{n, \mathfrak{X}_{K}/Y}$ 
under the map (\ref{MM}).
We define the subsheaf $V^{q}M^{n}_{r}$ ($q\geq 0$)
of $M^{n}_{r}$ as the part which is generated by 
$U^{q+1}M^{n}_{r}$
and
the image of 
$U^{p}\mathcal{K}^{M}_{n-1, \mathfrak{X}_{K}/Y}\otimes\langle\pi\rangle$
under the map $(\ref{MM})$.
\end{itemize}
\end{Definition}
Then we give a brief review 
of the structure of $M^{n}_{r}$.
Let $\omega_{Y}^{*}$ be 
the modified differential modules which is defined in
\cite[p.546, (1.5)]{H} (See also \cite[p.531]{SaP}).
For $q\geq 0$,
\begin{equation*}
\operatorname{gr}_{U/V}^{q}M^{n}_{r}
:=U^{q}M^{n}_{r}/V^{q}M^{n}_{r}
~~\textrm{and}~~
\operatorname{gr}_{V/U}^{q}M^{n}_{r}
:=V^{q}M^{n}_{r}/U^{q+1}M^{n}_{r}
\end{equation*}
are expressed by using
subsheaves of the modified differential modules $\omega_{Y}^{*}$
as follows:
\begin{Theorem}\upshape(Bloch-Kato \cite[pp.112--113, Corollary (1.4.1)]{B-K}
/Hyodo \cite[p.548, (1.7) Corollary]{H},
Sato \cite[pp.184--185, Theorem 3.3]{SaR}
)\label{QH}

Let the notations be the same as above. Then
\begin{itemize}
\item[(1)] The map (\ref{MM}) 
is surjective, 
that is,
the subsheaf 
$U^{0}M^{n}_{r}$
is the full sheaf 
$M^{n}_{r}$
for any $n\geq 0$ and $r>0$.
\item[(2)] Let $e$ be the absolute ramification index
of $K$, and let $r=1$. Then for $q$ with
$1\leq q< e^{\prime}:=pe/(p-1)$,
there are isomorphisms
\begin{align*}
\operatorname{gr}^{q}_{U/V}
M^{n}_{1}
\simeq
\left\{
\begin{array}{l}
\omega^{n-1}_{Y}/\mathcal{B}^{n-1}_{Y} ~~ ((p, q)=1),
\\
\omega^{n-1}_{Y}/\mathcal{Z}^{n-1}_{Y} ~~ (p|q),
\end{array}
\right.
\end{align*} 
\begin{align*}
\operatorname{gr}^{q}_{V/U}M^{n}_{1}
\xrightarrow{\sim}
\omega^{n-2}_{Y}/\mathcal{Z}^{n-2}_{Y},
\end{align*}
where
\begin{math}
\omega^{m}_{Y}
=W_{1}\omega^{m}_{Y}
\end{math}
and 
$\mathcal{B}^{m}_{Y}$ (resp.$\mathcal{Z}^{m}_{Y}$)
denotes the image of
$d: \omega_{Y}^{m-1}\to \omega_{Y}^{m}$
(resp. the kernel of 
$d: \omega_{Y}^{m}\to\omega_{Y}^{m+1}$).
\item[(3)] We have
\begin{equation*}
U^{q}M^{n}_{1}
=
V^{q}M^{n}_{1}
=0
\end{equation*}
for any $q\geq e^{\prime}$.
\end{itemize}
\end{Theorem}
We define the \'{e}tale subsheaf $FM^{n}_{r}$
as the part which is generated by
$U^{1}M^{n}_{r}$ and the image of
$(\iota^{*}\mathcal{O}^{\times}_{\mathfrak{X}})^{\otimes n}$ under the map (\ref{MM}).
\begin{Theorem}\upshape
(Sato \cite[p.533, Theorem 3.4.2]{SaP}, 
\cite[p.186, Theorem 3.4]{SaR}
)\label{QS}
Let $\mathfrak{X}$ be the same as in Theorem \ref{QH}.
Then there exists a short exact sequence of sheaves on 
$Y_{\et}$
\begin{equation*}
0
\to
FM^{n}_{r}
\to
M^{n}_{r}
\xrightarrow{\sigma^{n}_{\mathfrak{X}, r}}
\nu^{n-1}_{Y, r}
\to 0,
\end{equation*}
where $\sigma^{n}_{\mathfrak{X}, r}$ is induced by the boundary
map of Galois cohomology groups 
(cf.\cite[p.530, (3.2.5)]{SaP}).
%
%
Furthermore there is an isomorphism
\begin{equation*}
FM^{n}_{r}/U^{1}M^{n}_{r}
\xrightarrow{\sim}
\lambda^{n}_{Y, r}
\end{equation*}
sending a symbol
\begin{equation*}
\{x_{1}, x_{2}, \cdots, x_{n}\} ~~ (x_{i}\in \iota^{*}\mathcal{O}_{\mathfrak{X}}^{\times})
\end{equation*}
to
\begin{equation*}
d\log(
\bar{x_{1}}\otimes\bar{x_{2}}\otimes \cdots \otimes \bar{x_{n}}
).
\end{equation*}
Here for a section
$x\in \iota^{*}\mathcal{O}_{\mathfrak{X}}^{\times}$,
$\bar{x}$ denotes its residue class in 
$\mathcal{O}_{Y}^{\times}$.
\end{Theorem}
\begin{Remark}\upshape\label{filtiT}
By Definition \ref{DefpTt} and Theorem \ref{QS}, we have an isomorphism
\begin{equation*}
\mathcal{H}^{n}(\iota^{*}\mathfrak{T}_{r}(n))
\simeq
FM^{n}_{r}
\end{equation*}
for a positive integer $r$. Hence
there is a finite filtration of subsheaves of
$\mathcal{H}^{n}(\iota^{*}\mathfrak{T}_{1}(n))$
\begin{align*}
\mathcal{H}^{n}(\iota^{*}\mathfrak{T}_{1}(n))
\supset
U^{1}M^{n}_{1}
\supset\cdots\supset
U^{q}M^{n}_{1}
\supset
V^{q}M^{n}_{1}
\supset
U^{q+1}M^{n}_{1}
\supset
\cdots
\end{align*}
where we have isomorphisms
\begin{align}\label{iTU1}
\mathcal{H}^{n}(
\iota^{*}\mathfrak{T}_{1}(n))/
U^{1}M_{1}^{n}
\simeq
\lambda_{Y, 1}^{n},
\end{align}
\begin{align*}
U^{q}M^{n}_{1}/V^{q}M^{n}_{1}
\simeq
\left\{
\begin{array}{l}
\omega^{n-1}_{Y}/\mathcal{B}^{n-1}_{Y} ~~ ((p, q)=1),
\\
\omega^{n-1}_{Y}/\mathcal{Z}^{n-1}_{Y} ~~ (p|q)
\end{array}
\right.
\end{align*} 
for $1\leq q< e^{\prime}:=pe/(p-1)$,
\begin{align*}
V^{q}M^{n}_{1}/U^{q+1}M^{n}_{1}
\xrightarrow{\sim}
\omega^{n-2}_{Y}/\mathcal{Z}^{n-2}_{Y}
\end{align*}
for $1\leq q< e^{\prime}$ and
\begin{equation*}
U^{q}M^{n}_{1}
=
V^{q}M^{n}_{1}
=0
\end{equation*}
for any $q\geq e^{\prime}$
by Theorem \ref{QH} and Theorem \ref{QS}.
Here $e$ is the absolute ramification index of $K$.
\end{Remark}
\subsection{Rigidity}
The following lemma is useful to compute \'{e}tale cohomology groups with values in $U^{1}M_{1}^{n}$.
\begin{Lemma}\upshape\label{VBZ}
Let $A$ be the henselization of a local ring of
a normal crossing variety over 
a field of characteristic $p>0$. 
Then
\begin{equation*}
\HO^{i}_{\et}(A, \mathcal{B}^{n}_{A})
=
\HO^{i}_{\et}(A, \mathcal{Z}^{n}_{A})
=0
\end{equation*}
for $i\geq 1$. Here 
$\mathcal{B}^{n}_{A}$ (resp.$\mathcal{Z}^{n}_{A}$)
is the image of
$d: \omega^{n-1}_{A}\to \omega^{n}_{A}$
(resp.the kernel of
$d: \omega^{n}_{A}\to \omega^{n+1}_{A}$).
\end{Lemma}
\begin{proof}\upshape
$\mathcal{B}^{n}_{A}$ and $\mathcal{Z}^{n}_{A}$
are locally free $\mathcal{O}_{\spec (A)}$-modules
after twisting with Frobenius.
Therefore the statement follows from \cite[p.114, III, Remark 3.8]{M}
and \cite[p.103, III, Lemma 2.15]{M}.
\end{proof}

Then we are able to prove the following main result of this paper by
applying Theorem \ref{LL}.

\begin{Theorem}\upshape\label{MR}
Let $B$ be a discrete valuation ring of mixed characteristic $(0, p)$, $R$ the henselization of a local ring of a semistable family over $\spec(B)$
and $k$ the residue field of $R$. Then we have an isomorphism
\begin{equation*}
\HO_{\et}^{n+1}(R, \mathfrak{T}_{r}(n))   
\simeq
\HO^{1}_{\et}(k, \lambda_{k, r}^{n})
\end{equation*}
where
$\mathfrak{T}_{r}(n)$ is the $p$-adic Tate twist.
\end{Theorem}
\begin{proof}\upshape
Let 
\begin{math}
\iota: \spec (A) \hookrightarrow
\spec (R)
\end{math}
be the special fiber.
Then we have an isomorphism 
\begin{equation}\label{isoA}
\HO^{n+1}_{\et}(R, \mathfrak{T}_{r}(n))
\simeq
\HO^{n+1}_{\et}(A, \iota^{*}\mathfrak{T}_{r}(n))
\end{equation}
by \cite[p.777, The proof of Proposition 2.2.b)]{Ge}.
Since $A$ has 
$p$-cohomological dimension at most $1$,
we have
\begin{equation*}
\HO^{s}_{\et} 
(
A,
\mathcal{H}^{t}(
\iota^{*}\mathfrak{T}_{r}(n)
)
)
=0
\end{equation*}
for $s\geq 2$ and any positive integer $t$.
Moreover, we have
\begin{equation*}
\mathcal{H}^{t}(
\iota^{*}\mathfrak{T}_{r}(n)
)   
=0
\end{equation*}
for $t\geq n+1$ by the definition of $\mathfrak{T}_{r}(n)$.
So we have an isomorphism
\begin{equation}\label{n1T}
\HO^{n+1}_{\et}(A, \iota^{*}\mathfrak{T}_{r}(n))
\simeq
\HO^{1}_{\et}(A,
\mathcal{H}^{n}(\iota^{*}\mathfrak{T}_{r}(n))
)
\end{equation}
by the spectral sequence
\begin{equation*}
E^{s, t}_{2}
=
\HO^{s}_{\et}
(
A,
\mathcal{H}^{t}
(\iota^{*}
\mathfrak{T}_{r}(n)
)
)   
\Rightarrow
E^{s+t}
=
\HO^{s+t}_{\et} 
(
A,
\iota^{*}
\mathfrak{T}_{r}(n)
).
\end{equation*}

By (\ref{isoA}) and (\ref{n1T}), we have an isomorphism
\begin{equation}\label{isohi}
\HO^{n+1}_{\et}(R, \mathfrak{T}_{r}(n))  
\simeq
\HO^{1}_{\et}(A, \mathcal{H}^{n}(\iota^{*}\mathfrak{T}_{r}(n)))
\end{equation}
for any positive integer $r$ and 
so it suffices to show that we have an isomorphism
\begin{equation}\label{TL}
\HO_{\et}^{1}\left(
A, \mathcal{H}^{n}(\iota^{*}\mathfrak{T}_{r}(n))
\right)   
\simeq
\HO^{1}_{\et}(A, \lambda_{A, r}^{n}).   
\end{equation}
Since we have
\begin{equation*}
\HO^{i}_{\et}(
A,
U^{1}M^{n}_{1})
=0
\end{equation*}
for $i\geq 1$ by Theorem \ref{QH} and Lemma \ref{VBZ},
we have an isomorphism
\begin{equation*}
\HO^{1}_{\et}(A, 
\mathcal{H}^{n}(\iota^{*}\mathfrak{T}_{1}(n))
)
\simeq
\HO^{1}_{\et}
(
A. \lambda^{n}_{A, 1}
)
\end{equation*}
by (\ref{iTU1}).
Therefore we have
an isomorphism
\begin{equation*}
\HO^{n+1}_{\et}
(
A,
\iota^{*}\mathfrak{T}_{1}(n)
)
\simeq
\HO^{1}_{\et}
(
A,
\lambda_{A, 1}^{n}
)
\end{equation*}
by (\ref{n1T}).
Since
\begin{equation*}
\mathcal{H}^{t}(\mathfrak{T}_{r}(n))=0
\end{equation*}
for $t>n$, we have
\begin{equation*}
\HO^{n+2}_{\et}
(
A,
\iota^{*}
\mathfrak{T}_{r}(n)
)
=0
\end{equation*}
by the similar argument as in the proof of the isomorphism (\ref{n1T}).
Hence the top row in
a commutative diagram
\footnotesize
\begin{equation}\label{CTL}
\xymatrix
{
&
\HO^{n+1}_{\et} 
(A,
\iota^{*}\mathfrak{T}_{r}(n)
)
\ar[r]\ar[d]
&
\HO^{n+1}_{\et} 
(A,
\iota^{*}\mathfrak{T}_{r+1}(n)
)
\ar[r]\ar[d]
&
\HO^{n+1}_{\et} 
(A,
\iota^{*}\mathfrak{T}_{1}(n)
)
\ar[d]\ar[r]
&
0
\\
0\ar[r]
&
\HO^{1}_{\et}
(
A,
\lambda^{n}_{A, r}
)
\ar[r]
&
\HO^{1}_{\et}
(
A,
\lambda^{n}_{A, r+1}
)
\ar[r]
&
\HO^{1}_{\et}
(
A,
\lambda^{n}_{A, 1}
)
&
}    
\end{equation}
\normalsize
is exact. Moreover the bottom row in the
commutative diagram (\ref{CTL}) is exact 
by Theorem \ref{LL}.
If the homomorphism (\ref{TL}) is an isomorphism in the case where $r\leq j$,
then 
the homomorphism (\ref{TL}) is an isomorphism in the case where $r=j+1$
by applying the snake lemma to the commutative
diagram (\ref{CTL}).
Therefore the statement follows by induction on $r$.
\end{proof}
\subsection{Gersten-type conjecture}
\begin{Proposition}\upshape\label{spsup}
Let $A$ be a discrete valuation ring of mixed characteristic $(0, p)$
and $X$ a semistable family over $\spec(A)$.
Let $R$ be the local ring at a point of $X$ and
$Z$ the closed fiber of $\spec(R)$.
Then we have a
spectral sequence
\begin{equation*}
E_{1}^{s, t}=
\bigoplus_{x\in \spec(R)^{(s)}\cap Z}
\HO_{x}^{s+t}(R_{\et}, \mathfrak{T}_{r}(n))
\Rightarrow
E^{s+t}=
\HO_{Z}^{s+t}(R_{\et}, \mathfrak{T}_{r}(n))
\end{equation*}
for any integer $n$.
\end{Proposition}
\begin{proof}\upshape
For a bounded below complex
$F^{\bullet}$ of \'{e}tale sheaves on $\spec(R)$,
we have a spectral sequence
\begin{equation}\label{spBO}
\bigoplus_{x\in \spec(R)^{(s)}}  
\HO_{x}^{s+t}(R_{\et}, F^{\bullet})
\Rightarrow
E^{s+t}=
\HO^{s+t}(R_{\et}, F^{\bullet})
\end{equation}
(cf.\cite[Part 1, \S 1]{C-H-K}). 
Let $i: Z\to\spec(R)$ be the inclusion of
the closed fiber of $\spec(R)$ and
\begin{math}
F^{\bullet}=i_{*}Ri^{!}\mathfrak{T}_{r}(n).
\end{math}
Then we have
\begin{equation}\label{Supi}
\HO_{x}^{j}(R, F^{\bullet})
=
\left\{
\begin{array}{lcr}
\HO_{x}^{j}(R, \mathfrak{T}_{r}(n)) & (\textrm{if}~~x\in Z)
\\
0 & (\textrm{if}~~x\notin Z)
\end{array}
\right.
\end{equation}
for any integer $j>0$.
Hence the statement follows
from the spectral sequence (\ref{spBO}) and (\ref{Supi}).
\end{proof}
\begin{Lemma}\upshape\label{Iml}
Let $n$ and $N$ be positive integers.
Let $i$ be an integer where $i=0$ or $i=1$.
Consider a spectral sequence
\begin{math}
E^{s, t}_{1}
\Rightarrow
E^{s+t}
\end{math}
which satisfies the following conditions:
\begin{itemize}
\item[(a)] 
If $s<0$ or $s>N$, 
then 
\begin{math}
E^{s, t}_{1}=0;    
\end{math}
\item[(b)] If $s+t=n$ or $s+t=n+1$, 
then
\begin{math}
E^{s, t}_{2}=0
\end{math}
for $s\neq i$.
\end{itemize}
Then we have
\begin{equation*}
E^{n}=E^{i, n-i}_{2}.    
\end{equation*}
\end{Lemma}
\begin{proof}\upshape
By the definition of the spectral sequence,
we have a filtration
\begin{equation*}
E^{n}=E^{n}_{0}
\supset
E^{n}_{1}
\supset\cdots
\supset
E^{n}_{N}
\supset
0
\end{equation*}
which satisfies
\begin{equation*}
E^{n}_{k}/E^{n}_{k+1}
\simeq 
E^{k, n-k}_{\infty}
\end{equation*}
for any integer $k$. So we have
\begin{equation*}
E^{n}_{i+1}=0 
~~\textrm{and}~~
E^{n}_{i}=E^{n}
\end{equation*}
by the conditions.
Moreover, we have 
\begin{equation*}
E^{i, n-i}_{2}=E^{i, n-i}_{\infty}    
\end{equation*}
by the conditions. Hence the statement follows.
\end{proof}
\begin{Theorem}\upshape\label{SupInj}
Let $A$ be a discrete valuation ring
of mixed characteristic $(0, p)$, 
$X$ a semistable family over $\spec(A)$
and $R=\mathcal{O}_{X, x}$ the local ring at a point $x$ of $X$.
Let $Z$ be the closed fiber of $\spec(R)$. Then 
the sequence
\begin{equation*}
0\to  
\HO_{Z}^{n+2}(R_{\et}, \mathfrak{T}_{r}(n))
\to
\bigoplus_{x\in Z^{(0)}}
\HO_{x}^{n+2}(R_{\et}, \mathfrak{T}_{r}(n))
\to
\bigoplus_{x\in Z^{(1)}}
\HO_{x}^{n+3}(R_{\et}, \mathfrak{T}_{r}(n))
\end{equation*}
is exact for any integer $r>0$.
\end{Theorem}
\begin{proof}\upshape
By \cite[p.540, Theorem 4.4.7]{SaP} and Corollary \ref{GeZ}, the sequence
\begin{equation*}
\bigoplus_{x\in Z^{(0)}}
\HO_{x}^{s}(R_{\et}, \mathfrak{T}_{r}(n))    
\to
\bigoplus_{x\in Z^{(1)}}
\HO_{x}^{s+1}(R_{\et}, \mathfrak{T}_{r}(n))  
\to
\cdots
\end{equation*}
is exact for $s\leq n+1$.
So the statement follows from  
Proposition \ref{spsup}
and Lemma \ref{Iml}.
\end{proof}
\begin{Remark}\upshape
Let $A$ be a discrete valuation ring of mixed characteristic $(0, p)$
and $R$ a local ring 
of a smooth scheme over $\spec(A)$.
Let $r$ be a positive integer.
Then the sequence
\begin{equation*}
\bigoplus_{x\in \spec(R)^{(0)}}
\HO_{x}^{s}(R_{\et}, \mathbb{Z}/p^{r}(n))    
\to
\bigoplus_{x\in \spec(R)^{(1)}}
\HO_{x}^{s+1}(R_{\et}, \mathbb{Z}/p^{r}(n))  
\to
\cdots
\end{equation*}
is exact for $s\leq n$ by \cite[p.774, Theorem 1.2]{Ge} and \cite{V2}.
Hence we can also give an another proof of \cite[Theorem 1.1]{SM}
by Lemma \ref{Iml}
and the spectral sequence (\ref{spBO}).
\end{Remark}
\section{Applications}
Throughout this section, $n$ is a non negative integer
and $r$ is a positive integer.
\begin{Lemma}\upshape\label{ID}
Let $R$ be a discrete valuation ring
of mixed characteristic $(0, p)$.
Then the homomorphism
\begin{equation}\label{CI}
\HO^{n+1}_{\et}
(
R,
\mathfrak{T}_{r}(n)
)
\to
\HO^{n+1}_{\et}
(
k(R),
\mu_{p^{r}}^{\otimes n}
)
\end{equation}
is injective.
\end{Lemma}
\begin{proof}\upshape
Let $\mathfrak{m}$ be the residue field of $R$.
We have a commutative diagram
\begin{equation*}
\xymatrix{
K^{M}_{n}(k(R))/p^{r}
\ar[r]\ar[d]
&
K^{M}_{n-1}(\kappa(\mathfrak{m}))/p^{r}
\ar[d]
\\
\HO^{n}_{\et}(
k(R), \mu_{p^{r}}^{\otimes n}
)\ar[r]
&
\HO^{0}_{\et}(
\kappa(\mathfrak{m}),
\nu^{n-1}_{\kappa(\mathfrak{m}), r}
)
}    
\end{equation*}
where $K^{M}_{n}(F)$ is the $n$-th Milnor $K$-group
for a field $F$.
Then the top horizontal arrow is 
induced by
the boundary map (cf.\cite[p.529, \S 3.1]{SaP})
and so is surjective.
Moreover the right vertical map is an isomorphism
by \cite[p.113, Theorem (2.1)]{B-K}.
Hence the homomorphism
\begin{equation*}
\HO^{n}_{\et}(
k(R),
\mu_{p^{r}}^{\otimes n}
)
\to
\HO^{0}_{\et}(
\kappa(\mathfrak{m}),
\nu^{n-1}_{\kappa(\mathfrak{m}), r}
)
\end{equation*}
is surjective and the homomorphism
\begin{equation*}
\HO^{n+1}_{\et}(
R,
\mathfrak{T}_{r}(n)
)
\to
\HO^{n+1}_{\et}(
R,
\tau_{\leq n}Rj_{*}\mu^{\otimes n}_{p^{r}}
)
\end{equation*}
is injective
where 
\begin{math}
j: \spec (k(R))\to \spec (R)
\end{math}
is the generic fiber. 
Since we have a distinguished triangle
\begin{equation*}
\cdots\to
\tau_{\leq n}Rj_{*}\mu_{p^{r}}^{\otimes n}
\to
\tau_{\leq n+1}Rj_{*}\mu_{p^{r}}^{\otimes n}
\to
R^{n+1}j_{*}\mu_{p^{r}}^{\otimes n}[-(n+1)]
\to\cdots
\end{equation*}
and
\begin{equation*}
\HO^{n+1}_{\et}(
R,
\tau_{\leq n+1}Rj_{*}\mu^{\otimes n}_{p^{r}}
)
=
\HO^{n+1}_{\et}(
k(R),
\mu_{p^{r}}^{\otimes n}
),    
\end{equation*}
the homomorphism
\begin{equation*}
\HO^{n+1}_{\et}(
R,
\tau_{\leq n}Rj_{*}\mu^{\otimes n}_{p^{r}}
)
\to
\HO^{n+1}_{\et}(
k(R),
\mu_{p^{r}}^{\otimes n}
)
\end{equation*}
is injective. Therefore the statement follows.
\end{proof}
\begin{Proposition}\upshape\label{GHI}
Let 
$\mathfrak{X}$ be a semistable family over
the spectrum of a discrete valuation ring
of mixed characteristic $(0, p)$
and
$R=\tilde{\mathcal{O}_{\mathfrak{X}, x}}$ 
be the henselian local ring of $x$ in $\mathfrak{X}$.
Then the homomorphism
(\ref{CI})
is injective.
\end{Proposition}
\begin{proof}\upshape
By \cite[Corollary 2.1.3]{C-H-K},
it suffices to prove the statement in the case where $x$ is in the closed fiber
of $\mathfrak{X}$.
Let 
$\mathfrak{m}$
be the maximal ideal of $R$. Let
$a\in\mathfrak{m}\setminus\mathfrak{m}^{2}$ 
and 
$p\in (a)=\mathfrak{p}$. 
By the assumption of $\mathfrak{X}$, 
there exists such a prime ideal $\mathfrak{p}=(a)$. 
Then 
$R/\mathfrak{p}$ is a regular local ring of characteristic $p>0$
and
we have a commutative diagram
\begin{equation}\label{CTN}
\xymatrix{
\HO^{n+1}_{\et}(
R, \mathfrak{T}_{r}(n)
)
\ar[r]\ar[d]
&
\HO^{n+1}_{\et}(
R_{\mathfrak{p}}, \mathfrak{T}_{r}(n)
) \ar[d] \\
\HO^{n+1}_{\et}(
R/\mathfrak{p},
\lambda_{R/\mathfrak{p}, r}^{n}
)
\ar[r]
&
\HO^{n+1}_{\et}(
\kappa(\mathfrak{p}),
\lambda_{\kappa(\mathfrak{p}), r}^{n}
)
}
\end{equation}
where the left vertical arrow in (\ref{CTN}) is injective by 
Theorem \ref{MR}.
Since $R/\mathfrak{p}$ is a regular local ring,
the bottom horizontal arrow in (\ref{CTN}) is also injective by
Theorem \ref{ninj}.
Hence the top horizontal arrow 
in (\ref{CTN}) is injective.
Therefore the statement follows from Lemma \ref{ID}.
\end{proof}
\begin{Lemma}\upshape\label{SI}
Let $R$ be a local ring of
a 
semistable family
over a discrete valuation ring 
of mixed characteristic $(0, p)$.
Suppose that 
$\operatorname{dim}(R)=2$.
Let $\mathfrak{p}$ be a prime ideal
such that $\operatorname{ht}(\mathfrak{p})=1$,
$R/\mathfrak{p}$ is regular and
$(p)\subset \mathfrak{p}$.
Then the homomorphism
\begin{equation*}
\HO^{n+2}_{R/\mathfrak{p}}
(
R,
\mathfrak{T}_{r}(n)
)
\to
\HO^{n+2}_{\mathfrak{p}}
(
R_{\mathfrak{p}},
\mathfrak{T}_{r}(n)
)
\end{equation*}
is injective.
\end{Lemma}
\begin{proof}\upshape
Let $\mathfrak{m}$ be the maximal ideal of $R$.
Since the sequence
\begin{equation*}
\HO^{n+1}_{\mathfrak{p}}
(
R_{\mathfrak{p}},
\mathfrak{T}_{r}(n)
)
\to 
\HO^{n+2}_{R/\mathfrak{m}}
(
R,
\mathfrak{T}_{r}(n)
) 
\to
\HO^{n+2}_{R/\mathfrak{p}}
(
R,
\mathfrak{T}_{r}(n)
)
\to
\HO^{n+2}_{\mathfrak{p}}
(
R_{\mathfrak{p}},
\mathfrak{T}_{r}(n)
)
\end{equation*}
is exact by \cite[p.92, III, Remark 1.26]{M}, it suffices to show that 
the homomorphism 
\begin{equation}\label{SMP}
\HO^{n+1}_{\mathfrak{p}}
(
R_{\mathfrak{p}},
\mathfrak{T}_{r}(n)
)
\to 
\HO^{n+2}_{R/\mathfrak{m}}
(
R,
\mathfrak{T}_{r}(n)
)
\end{equation}
is surjective.
The homomorphism (\ref{SMP}) coincides with
$\pm$ the boundary map
\begin{equation}\label{bH}
\HO^{0}_{\et}(
k(R/\mathfrak{p}),
\nu^{n-1}_{r}
)
\to
\HO^{0}_{\et}(
R/\mathfrak{m},
\nu^{n-2}_{r}
)
\end{equation}
by \cite[p.547, Theorem 6.1.1]{SaP}, \cite[p.540, Theorem 4.4.7]{SaP} and 
\cite[p.187, Remark 3.7]{SaR}.
Let
$K^{M}_{n}(F)$ be
the $n$-th Milnor $K$-group for a field $F$.
Then the homomorphism (\ref{bH}) coincides with
$\pm$ the tame symbol
\begin{equation*}
K^{M}_{n-1}(k(R/\mathfrak{p}))/p^{r}
\to
K^{M}_{n-2}(R/\mathfrak{m})/p^{r}
\end{equation*}
of Milnor $K$-groups (modulo $p^{r}$)
by the definition of the boundary map (cf.\cite[p.150, (1.3) (ii)]{K}) and
\cite[p.113, Theorem (2.1)]{B-K}. So the homomorphism (\ref{SMP}) is surjective. 
This completes the proof.
\end{proof}
\begin{Lemma}\upshape\label{BLNe}
Let $B$ be a discrete valuation ring of mixed 
characteristic $(0, p)$. 
Let $\mathfrak{X}$ be a semistable family over $\spec(B)$.
Then we have 
an isomorphism in 
$D^{b}(\mathfrak{X}_{\et}, \mathbb{Z}/p^{r}\mathbb{Z})$
\begin{equation*}
\tau_{\leq n} \left(
\mathbb{Z}/p^{r}(n)_{\Nis}
\right)
\simeq
\tau_{\leq n} \left(
R\alpha_{*}
\mathbb{Z}/p^{r}(n)_{\et}
\right),
\end{equation*}
where $\alpha: \mathfrak{X}_{\et}\to \mathfrak{X}_{\Nis}$ is 
the change of site map
and $D^{b}(\mathfrak{X}_{\et}, \mathbb{Z}/p^{r}\mathbb{Z})$
is the derived category of bounded complexes of \'{e}tale 
$\mathbb{Z}/p^{r}\mathbb{Z}$-sheaves on $\mathfrak{X}$.
\end{Lemma}
\begin{proof}\upshape
Let $\epsilon: \mathfrak{X}_{\et}\to \mathfrak{X}_{\Zar}$
and
$\beta: \mathfrak{X}_{\Nis}\to \mathfrak{X}_{\Zar}$
be the change of site map.  Let 
$j:U\to \mathfrak{X}$ 
be
the generic fiber
and 
$i:Z\to \mathfrak{X}$
be the closed fiber.
Then the Beilinson-Lichtenbaum conjecture holds for $U$ 
by \cite[Theorem 1.1]{G-L2} and the Bloch-Kato conjecture (\cite{V2}). 
Moreover, 
the Beilinson-Lichenbaum conjecture holds for $Z$
by Proposition \ref{apbk}.2. 
Therefore we have a quasi-isomorphism
\begin{equation}\label{ZeZ}
\tau_{\leq n} \left(
\mathbb{Z}/p^{r}(n)_{\Zar}
\right)
\simeq
\tau_{\leq n} \left(
R\epsilon_{*}
\mathbb{Z}/p^{r}(n)_{\et}
\right)
\end{equation}
by 
\cite[p.33, Proposition 2.1]{SM} and
the same argument as in the proof of 
\cite[Theorem 1.2.2]{Ge}.
Moreover we have quasi-isomorphisms
\begin{equation*}
\beta^{*}R\beta_{*}
\simeq
\operatorname{id}
~~
\textrm{and}
~~
R\epsilon_{*}
\simeq
R\beta_{*}R\alpha_{*}
.
\end{equation*}
Since $\beta^{*}$ is exact and 
\begin{math}
\beta^{*}(\mathbb{Z}(n)_{\Zar})
\simeq
\mathbb{Z}(n)_{\Nis},
\end{math}
the statement follows from (\ref{ZeZ}).
\end{proof}
\begin{Proposition}\upshape\label{VAM}
Let $B$ be a discrete valuation ring of mixed 
characteristic $(0, p)$ and
$R$ the henselization of a local ring of a semistable family over
$\spec(B)$.
Then we have
\begin{equation*}
\HO^{n+1}_{\Zar}(
R, 
\mathbb{Z}/p^{r}(n)
)
=0
\end{equation*}
for any integer $r>0$
where $\mathbb{Z}(n)$ is 
Bloch's cycle complex
for the Zariski topology
(\cite{B}, \cite{L})
and 
\begin{math}
\mathbb{Z}/p^{r}(n)
=
\mathbb{Z}(n)\otimes\mathbb{Z}/p^{r}.
\end{math}
\end{Proposition}
\begin{proof}\upshape
Let $Z$ be the closed fiber of $\spec(R)$
and $U$ the generic fiber
of $\spec(R)$.
By Proposition \ref{GHI},
the homomorphism
\begin{equation*}
\HO^{n+1}_{\et}(R, \mathfrak{T}_{r}(n))    
\to
\HO^{n+1}_{\et}(U, \mu_{p^{r}}^{\otimes n})
\end{equation*}
is injective
and so the homomorphism
\begin{equation*}
\HO^{n}_{\et}(U, \mu_{p^{r}}^{\otimes n})
\to
\HO^{n+1}_{\et}(Z, \nu_{r}^{n-1}[-(n+1)])
\end{equation*}
is surjective by \cite[p.540, Theorem 4.4.7]{SaP}.
So the homomorphism
\begin{equation*}
\HO^{n}_{\Zar}(U, \mathbb{Z}/p^{r}(n))
\to
\HO^{n-1}_{\Zar}(Z, \mathbb{Z}/p^{r}(n-1))
\end{equation*}
is surjective by 
\cite[p.774, Theorem 1.2.2]{Ge}, \cite{V2} and Proposition \ref{iso}. 
Hence the homomorphism
\begin{equation}\label{ZUkin}
\HO_{\Zar}^{n+1}(R, \mathbb{Z}/p^{r}(n))
\to
\HO_{\Zar}^{n+1}(U, \mathbb{Z}/p^{r}(n))
\end{equation}
is injective by the localization theorem 
(\cite[p.780, Corollary 3.3 a)]{Ge}).
Moreover, we have the commutative diagram
\footnotesize
\begin{equation}\label{Ukmi}
\xymatrix{
&
\HO^{n+1}_{\et}(R,
\mathfrak{T}_{r}(n))
\ar[r]\ar[d]
&
\HO^{n+1}_{\et}(U, \mu_{p^{r}}^{\otimes n})
\ar[r]\ar[d] 
&
\HO_{Z}^{n+2}(R_{\et}, \mathfrak{T}_{r}(n))  \ar[d]
\\
0 \ar[r]
&
\displaystyle\bigoplus_{x\in Z^{(0)}}
\HO^{n+1}_{\et}(R_{x}, \mathfrak{T}_{r}(n))
\ar[r]
&
\displaystyle\bigoplus_{x\in Z^{(0)}}
\HO^{n+1}_{\et}(k(R), \mu_{p^{r}}^{\otimes n})
\ar[r]
&
\displaystyle\bigoplus_{x\in Z^{(0)}}
\HO_{x}^{n+2}(R_{\et}, \mathfrak{T}_{r}(n))
}    
\end{equation}
\normalsize
where the sequences are exact. 
Then the left map in the diagram (\ref{Ukmi}) is injective 
by Proposition \ref{GHI} and 
the right map in the diagram (\ref{Ukmi}) is injective by 
Theorem \ref{SupInj}.
So the middle map in the diagram (\ref{Ukmi}) is injective
and the homomorphism
\begin{equation*}
\HO^{n+1}_{\et}(U, \mathbb{Z}/p^{r}(n))
\to
\HO^{n+1}_{\et}(k(R), \mathbb{Z}/p^{r}(n))
\end{equation*}
is injective by \cite[p.774, Theorem 1.2.4]{Ge} and \cite{V2}. 
Moreover, the homomorphism
\begin{equation*}
\HO^{n+1}_{\Zar}(U, \mathbb{Z}/p^{r}(n)) 
\to
\HO^{n+1}_{\et}(U, \mathbb{Z}/p^{r}(n)) 
\end{equation*}
is injective by \cite[Theorem 1.6]{G-L2} and \cite{V2}.
Hence the homomorphism 
\begin{equation*}
\HO_{\Zar}^{n+1}(U, \mathbb{Z}/p^{r}(n))\to 
\HO_{\Zar}^{n+1}(k(R), \mathbb{Z}/p^{r}(n))
\end{equation*}
is injective. 
By the definition of $\mathbb{Z}/p^{r}(n)$
and the same argument in the proof of Proposition \ref{van}, 
we have
\begin{equation*}
\HO^{n+1}_{\Zar}(k(R), \mathbb{Z}/p^{r}(n)) 
=
\HO^{0}_{\Zar}(
k(R), \mathcal{H}^{n+1}(\mathbb{Z}/p^{r}(n))
)
=0
\end{equation*}
and so we have
\begin{equation}\label{ZvanU}
\HO_{\Zar}^{n+1}(U, \mathbb{Z}/p^{r}(n))=0.    
\end{equation}
%
Therefore the statement follows from the equation (\ref{ZvanU}) and 
the injectivity of the homomorphism (\ref{ZUkin}).
\end{proof}
The following means that 
\cite[p.524, Conjecture 1.4.1 (1)]{SaP} 
holds in the case where $\mathfrak{X}$ is a semistable family
over the spectrum of a Dedekind domain
and $\operatorname{dim}(\mathfrak{X})=2$.
\begin{Corollary}\upshape
Let $\mathfrak{X}$ be a semistable family
over the spectrum of a discrete valuation ring 
of mixed characteristic $(0, p)$.
Suppose that 
$\operatorname{dim}(\mathfrak{X})=2$.
Then we have an isomorphism in 
$D^{b}(\mathfrak{X}_{\et}, \mathbb{Z}/p^{r}\mathbb{Z})$
\begin{equation*}
\mathbb{Z}/p^{r}(n)_{\et}
\simeq
\mathfrak{T}_{r}(n),
\end{equation*}
where
$\mathbb{Z}/p^{r}(n)_{\et}$
is the \'{e}tale sheafification of
$\mathbb{Z}/p^{r}(n)$
and
$D^{b}(\mathfrak{X}_{\et}, \mathbb{Z}/p^{r}\mathbb{Z})$
is the derived category of bounded complexes of
\'{e}tale $\mathbb{Z}/p^{r}\mathbb{Z}$-sheaves on $\mathfrak{X}$.
\end{Corollary}
\begin{proof}\upshape
We have a quasi-isomorphism
\begin{equation*}
\tau_{\leq n}\left(
\mathbb{Z}/p^{r}(n)_{\et}
\right)
\simeq
\tau_{\leq n}
\left(
\mathfrak{T}_{r}(n)
\right)
\end{equation*}
by \cite[p.209, Remark 7.2]{SaR}.
By the definition of 
$\mathfrak{T}_{r}(n)$, 
we have
\begin{equation*}
\mathcal{H}^{i}(\mathfrak{T}_{r}(n))
=0
\end{equation*}
for $i\geq n+1$.
Let 
\begin{math}
\alpha: \mathfrak{X}_{\et}
\to \mathfrak{X}_{\Nis}
\end{math}
be the canonical map of sites.
Then the sheafification $\alpha^{*}$ is exact and
\begin{equation*}
\alpha^{*}\mathbb{Z}/p^{r}(n)_{\Nis}
=
\mathbb{Z}/p^{r}(n)_{\et}.    
\end{equation*}
Hence it suffices to show that
\begin{equation*}
\mathcal{H}^{i} 
(
\mathbb{Z}/p^{r}(n)
)
=0
\end{equation*}
for $i\geq n+1$. 
Since we have
\begin{align*}
\Gamma%
\left
(
\spec (\tilde{\mathcal{O}_{\mathfrak{X}, x}}),
\mathcal{H}^{n+1}(\mathbb{Z}/p^{r}(n)_{\Nis})
\right
)
&\simeq
\HO^{n+1}%
\left
(
\mathbb{Z}/p^{r}(n)_{\Zar}
(\spec(
\tilde{\mathcal{O}_{\mathfrak{X}, x}}
))
\right
) \\
&\simeq
\HO_{\Zar}^{n+1}
(
\tilde{
\mathcal{O}_{\mathfrak{X}, x}
},
\mathbb{Z}/p^{r}(n)
) \nonumber
\end{align*}
by \cite[p.779, Theorem 3.2 b)]{Ge},
it suffices to show that
\begin{equation*}
\HO^{i}_{\Zar}(
\tilde{
\mathcal{O}_{\mathfrak{X}, x}
},
\mathbb{Z}/p^{r}(n)
)
=0
\end{equation*}
for $i\geq n+1$.
Here $\tilde{\mathcal{O}_{\mathfrak{X}, x}}$
is the henselian local ring at a point $x$
of $\mathfrak{X}$.

Put
$R=\tilde{\mathcal{O}_{\mathfrak{X}, x}}$.
Let
$Z\hookrightarrow \spec(R)$ 
be a closed immersion of codimension $1$
with 
$U=\spec(R)\setminus Z$. Suppose that $Z$ is regular and $\operatorname{char}(Z)=p>0$.
By the assumptions of $\mathfrak{X}$, we are able to choose such a $Z$.
Then we have
\begin{equation*}
\HO^{i}_{\Zar}(
Z, 
\mathbb{Z}/p^{r}(n-1)
)
=0
\end{equation*}
for $i> n$ by 
\cite[p.786, Corollary 4.4]{Ge}.
Since $\operatorname{dim}(U)\leq1$, we have
\begin{equation*}
\HO^{i}\left(
\mathbb{Z}(n)_{\Zar}(U)
\right)
=0
\end{equation*}
for $i>n+1$
by the definition of $\mathbb{Z}(n)_{\Zar}$. Moreover we have
\begin{equation*}
\HO_{\Zar}^{i}\left(
U, \mathbb{Z}(n)
\right)
=
\HO^{i}
\left(
\mathbb{Z}(n)_{\Zar}(U)
\right)
=0    
\end{equation*}
for $i>n+1$ 
by \cite[p.779, Theorem 3.2 b)]{Ge}.
Hence we have
\begin{equation*}
\HO^{i}_{\Zar}(
U,
\mathbb{Z}/p^{r}(n)
)
=0
\end{equation*}
for $i>n+1$.
Therefore
we have
\begin{equation*}
\HO^{i}_{\Zar}(R,
\mathbb{Z}/p^{r}(n))
=0
\end{equation*}
for $i>n+1$ 
by the localization theorem (\cite[p.779, Theorem 3.2 a)]{Ge}). 
Moreover we have
\begin{equation*}
\HO^{n+1}_{\Zar}
(R,
\mathbb{Z}/p^{r}(n))
=0   
\end{equation*}
by Proposition \ref{VAM}. 
This completes the proof.
\end{proof}
\begin{Remark}\label{BGl}\upshape
Let $\mathfrak{X}$ be a regular scheme 
which is finite type over the spectrum of
a Dedekind domain.
Let $R$ be the henselization of 
a local ring $\mathcal{O}_{\mathfrak{X}, x}$
and 
$l$ a positive integer which is invertible in $R$.
Then we can prove
\begin{equation*}
\HO^{n+1}_{\Zar}(R, \mathbb{Z}/l(n))
=
\HO^{n+1}_{\Nis}(R, \mathbb{Z}/l(n))
=0    
\end{equation*}
by the similar argument as in Proposition \ref{VAM}. Hence we have
\begin{equation*}
\mathcal{H}^{n+1}(\mathbb{Z}/l(n)_{\et})
=
0
\end{equation*}
by \cite[p.779, Theorem 3.2.b)]{Ge}.
Let $\iota: Z\hookrightarrow \spec(R)$ be a closed immersion
of regular
subschemes of pure codimension $1$. Then we have an isomorphism
\begin{equation*}
R\iota^{!}\mu_{l}^{\otimes n} 
\simeq
\mu_{l}^{\otimes (n-1)}[-2]
\end{equation*}
in $D^{b}(Z_{\et}, \mathbb{Z}/l\mathbb{Z})$
by the absolute purity theorem (\cite{G}).
Here $D^{b}(Z_{\et}, \mathbb{Z}/l\mathbb{Z})$ 
is the derived category of bounded
complexs of \'{e}tale $\mathbb{Z}/l\mathbb{Z}$-sheaves 
on $Z$.
On the other hand, we have an isomorphism
\begin{equation*}
\tau_{\leq n+1}(R\iota^{!}\mathbb{Z}/l(n)_{\et}^{\spec(R)})
\simeq
\tau_{\leq n+1}(\mathbb{Z}/l(n-1)_{\et}^{Y}[-2])
\end{equation*}
in $D^{b}(Z_{\et}, \mathbb{Z}/l\mathbb{Z})$
by \cite[p.33, Proposition 2.1]{SM} and \cite{V2}.
Hence we have an isomorphism
\begin{equation*}
\tau_{\leq n}(\mathbb{Z}/l(n)_{\et})
\simeq
\mu_{l}^{\otimes n}
\end{equation*}
in $D^{b}(\spec(R)_{\et}, \mathbb{Z}/l\mathbb{Z})$
by induction on $\operatorname{dim}(R)$.
Therefore we have an isomorphism
\begin{equation*}
\tau_{\leq n+1} 
(\mathbb{Z}/l(n)_{\et})
\simeq
\mu_{l}^{\otimes n}
\end{equation*}
in $D^{b}(\mathfrak{X}_{\et}, \mathbb{Z}/l\mathbb{Z})$.
\end{Remark}
\begin{Theorem}\upshape\label{LGe}
Let 
$\mathfrak{X}$ be a semistable family
over the spectrum of
a discrete valuation ring 
of mixed characteristic $(0, p)$ and
$R=\tilde{\mathcal{O}_{\mathfrak{X}, x}}$ 
the henselian local ring
of $x$ in $\mathfrak{X}$. 
Suppose that $\operatorname{dim}(R)=2$.
Then the sequence
%
\begin{equation*}
0\to 
\HO^{n+1}_{\et}(
R,
\mathfrak{T}_{r}(n)
)
\to
\HO^{n+1}_{\et}
(
k(R),
\mu_{p^{r}}^{\otimes n}
)
\xrightarrow{(*)}
\bigoplus_{\substack{
\mathfrak{p}\in\spec (R)^{(1)}
\\
}
}
\HO^{n+2}_{\mathfrak{p}}
\left(
(R_{\mathfrak{p}})_{\et},
\mathfrak{T}_{r}(n)
\right)
\end{equation*}
%
is exact for any integers $n\geq 0$ and $r>0$.
\end{Theorem}
\begin{proof}\upshape
Let $\mathfrak{m}$ be the maximal ideal of $R$.
Let 
$a\in\mathfrak{m}\setminus\mathfrak{m}^{2}$ and $p\in (a)$.
By the assumption of $\mathfrak{X}$, there exists
such an element $a\in \mathfrak{m}\setminus\mathfrak{m}^{2}$.
Let us denote $\mathfrak{q}=(a)$ and
$U=\spec (R)\setminus \spec (R/\mathfrak{q})$. 
Put
\begin{equation*}
\HO^{n+2}_{R/\mathfrak{q}}(R_{\et}, \mathfrak{T}_{r}(n))^{\prime}
=
\operatorname{Im}\left(
\HO^{n+1}_{\et}(U, \mathfrak{T}_{r}(n))
\to
\HO^{n+2}_{R/\mathfrak{q}}(R_{\et}, \mathfrak{T}_{r}(n))
\right)
\end{equation*}
and
\begin{equation*}
\HO^{n+1}_{\et}(
k(R),
\mu^{\otimes n}_{p^{r}}
)^{\prime}   
=
\operatorname{Ker}
\left(
\HO^{n+1}_{\et}(
k(R),
\mu^{\otimes n}_{p^{r}}
)
\to
\bigoplus_{
\mathfrak{p}\in U^{(1)}
}
\HO^{n+2}_{\mathfrak{p}}
\left(
(R_{\mathfrak{p}})_{\et}, 
\mathfrak{T}_{r}(n)
\right)
\right).
\end{equation*}
Then we have a 
commutative diagram
\footnotesize
\begin{equation}
\xymatrix{
&
\HO^{n+1}_{\et}(
R, 
\mathfrak{T}_{r}(n)
)\ar[d]\ar[r]
&
\HO^{n+1}_{\et}(
U, 
\mathfrak{T}_{r}(n)
)\ar[r]\ar[d]
&
\HO^{n+2}_{R/\mathfrak{q}}
(
R_{\et}, 
\mathfrak{T}_{r}(n)
)^{\prime}
\ar[r]\ar[d]
&
0
\\
0
\ar[r]
&
\operatorname{Ker}
\left(
(*)
\right)
\ar[r]
&
\HO^{n+1}_{\et}(
k(R),
\mu^{\otimes n}_{p^{r}}
)^{\prime}\ar[r]
&
\HO^{n+2}_{\mathfrak{q}}
\left(
(R_{\mathfrak{q}})_{\et}, 
\mathfrak{T}_{r}(n)
\right)
}    
\end{equation}
\normalsize

where the rows are exact. Then the middle vertical arrow is surjective
by \cite[p.778, Lemma 2.4]{Ge}. 
Moreover the right vertical arrow is injective by Lemma \ref{SI}. 
Hence the surjectivity of the left vertical arrow follows from the snake lemma.
Moreover the left vertical arrow is injective by Proposition \ref{GHI}. This completes the proof.
\end{proof}
\begin{Corollary}\upshape\label{LGLC}
Let $R$ be the same as in Theorem \ref{LGe}
and $Y$ the special fiber of $\spec (R)$.
Then the natural map
%
\begin{equation*}
\HO^{n+1}_{\et}
(
k(R),
\mu_{p^{r}}^{\otimes n}
)
\to
\bigoplus_{\mathfrak{p}\in \spec (R)^{(1)}\setminus Y^{(0)}}
\HO^{n}_{\et}
(
\kappa(\mathfrak{p}),
\mu_{p^{r}}^{\otimes (n-1)}
)
\oplus
\bigoplus_{
\mathfrak{p}\in Y^{(0)}
}
\HO^{n+1}_{\et}
(
k(\tilde{R_{\mathfrak{p}}})
,
\mu_{p^{r}}^{\otimes n}
)
\end{equation*}
%
%
is injective for any integers $n\geq 0$ and $r>0$.
Here $\tilde{R_{\mathfrak{p}}}$ is the henselization of 
a local ring $R_{\mathfrak{p}}$.
\end{Corollary}
\begin{proof}\upshape
Put
\begin{equation*}
\operatorname{RH}^{n+1}(R)
:=
 \bigoplus_{\mathfrak{p}\in \spec (R)^{(1)}\setminus Y^{(0)}
}
\HO^{n}_{\et}
(
\kappa(\mathfrak{p}),
\mu_{p^{r}}^{\otimes (n-1)}
)
\oplus
\bigoplus_{
\mathfrak{p}\in Y^{(0)}
}
\HO^{n+1}_{\et}
(
k(\tilde{R_{\mathfrak{p}}})
,
\mu_{p^{r}}^{\otimes n}
).
\end{equation*}
Then we have a commutative diagram
\footnotesize
\begin{equation*}
\xymatrix{
0
\ar[r]
&\HO^{n+1}_{\et}(
R
,
\mathfrak{T}_{r}(n)
)
\ar[r]\ar[d]
&\HO^{n+1}_{\et}(
k(R)
,
\mathfrak{T}_{r}(n)
)
\ar[r]\ar[d]
&
\displaystyle\bigoplus_{\mathfrak{p}\in \spec(R)^{(1)}}
\HO^{n+2}_{\mathfrak{p}}\left(
(R_{\mathfrak{p}})_{\et},
\mathfrak{T}_{r}(n)   
\right)\ar@{=}[d]
\\
0
\ar[r]
&
\displaystyle
\bigoplus_{
\mathfrak{p}\in Y^{(0)}
}
\HO^{n+1}_{\et}(
\tilde{R_{\mathfrak{p}}},
\mathfrak{T}_{r}(n)
)
\ar[r]
&
\operatorname{RH}^{n+1}(R)
\ar[r]
&
\displaystyle\bigoplus_{\mathfrak{p}\in \spec(R)^{(1)}}
\HO^{n+2}_{\mathfrak{p}}\left(
(R_{\mathfrak{p}})_{\et},
\mathfrak{T}_{r}(n)   
\right),
}
\end{equation*}
\normalsize
where the rows are exact by Theorem \ref{LGe} and 
the absolute purity theorem
(cf.\cite[p.241, VI, Theorem 5.1]{M}).
Hence it suffices to show that
the homomorphism
\begin{equation*}
\HO^{n+1}_{\et}(
R, \mathfrak{T}_{r}(n)
)
\to
\bigoplus_{
\mathfrak{p}\in Y^{(0)}
}
\HO^{n+1}_{\et}(
\tilde{R_{\mathfrak{p}}},
\mathfrak{T}_{r}(n)
)
\end{equation*}
is injective.
By the assumption of $R$, there exists an element
$\mathfrak{p}\in Y^{(0)}$ such that 
$R/\mathfrak{p}$ is regular.
Moreover, a diagram
\begin{equation*}
\xymatrix{
\HO^{n+1}_{\et}(
R, 
\mathfrak{T}_{r}(n)
)
\ar[r]\ar[d]
&
\HO^{n+1}_{\et}(
\tilde{R_{\mathfrak{p}}}, 
\mathfrak{T}_{r}(n)
)
\ar[d]
\\
\HO^{1}_{\et}(
R/\mathfrak{p},
\lambda_{r}^{n}
)
\ar[r]
&
\HO^{1}_{\et}(
\kappa(\mathfrak{p}),
\lambda_{r}^{n}
)
}    
\end{equation*}
is commutative. Hence
the homomorphism
\begin{equation*}
\HO^{n+1}_{\et}(
R, \mathfrak{T}_{r}(n)
)
\to
\HO^{n+1}_{\et}(
\tilde{R_{\mathfrak{p}}},
\mathfrak{T}_{r}(n)
)
\end{equation*}
is injective 
by Theorem \ref{MR} and 
Theorem \ref{ninj}.
This completes the proof.
\end{proof}
\begin{Lemma}\upshape\label{WA}
Let $v_{i}$ $(1\leq i\leq m)$ be a finite
collection of independent discrete valuations on
a field $K$ of characteristic $0$.
Let $K_{i}$ be the henselization of $K$
at $v_{i}$ and
$\kappa(v_{i})$ the residue field
of $v_{i}$ for each $i$. 
Suppose that 
$\operatorname{char}(\kappa(v_{i}))=p>0$
for all $i$.
Then the natural map
\begin{equation*}
\HO^{n}_{\et}(
K, \mu_{p^{r}}^{\otimes n}
)
\to
\bigoplus_{i}
\HO^{n}_{\et}(
K_{i}, 
\mu_{p^{r}}^{\otimes n}
)
\end{equation*}
is surjective for integers $n\geq 1$ and $r>0$.
\end{Lemma}
\begin{proof}\upshape
The statement follows from \cite[p.131, Theorem (5.12)]{B-K} 
and \cite[pp.61--62, Lemma 6.2]{SM}.
\end{proof}
\begin{Theorem}\upshape\label{MTP}
Let $\mathfrak{X}$ be a proper and semistable
family over an excellent henselian discrete valuation ring of mixed
characteristic $(0, p)$
and $Y$ the special fiber of $\mathfrak{X}$. 
Let $r$ be any positive integer.
Suppose that $\operatorname{dim}(\mathfrak{X})=2$.
Then the map
\begin{equation*}
\HO^{n+1}_{\et}
(
k(\mathfrak{X}),
\mu_{p^{r}}^{\otimes n}
)
\to
\bigoplus_{\mathfrak{p}\in 
\mathfrak{X}^{(1)}\setminus Y^{(0)}
}
\HO^{n}_{\et}
(
\kappa(\mathfrak{p}),
\mu_{p^{r}}^{\otimes (n-1)}
)
\oplus
\bigoplus_{\mathfrak{p}\in Y^{(0)}
}
\HO^{n+1}_{\et}
(
k(\tilde{\mathcal{O_{\mathfrak{X}, \mathfrak{p}}}})
,
\mu_{p^{r}}^{\otimes n}
)
\end{equation*}
%
is injective for $n\geq 1$.
Here 
$\tilde{\mathcal{O}_{\mathfrak{X}, \mathfrak{p}}}$
is the henselization of a local ring
$\mathcal{O}_{\mathfrak{X}, \mathfrak{p}}$.
This implies that the local-global map
\begin{equation*}
\HO^{n+1}_{\et}
(
k(\mathfrak{X}),
\mu_{p^{r}}^{\otimes n}
)
\to
\prod_{\mathfrak{p}\in \mathfrak{X}^{(1)}}
\HO^{n+1}_{\et}
(
k(\tilde{\mathcal{O_{\mathfrak{X}, \mathfrak{p}}}})
,
\mu_{p^{r}}^{\otimes n}
)
\end{equation*}
%
is injective for $n\geq 1$.
\end{Theorem}
\begin{proof}
The proof of the statement is the same as \cite[Theorem 2.5]{Hu}.
The statement follows from Corollary \ref{LGLC} 
and Lemma \ref{WA}.
We review the proof of \cite[Theorem 2.5]{Hu}
for convenience. 

Let
\begin{equation*}
\xymatrix{
X_{1}\ar[d]_{j_{1}} 
&
X_{3}\ar[l]_{j^{\prime}_{2}}\ar[d]^{j^{\prime}_{1}} \\
X_{0} 
&
X_{2}\ar[l]^{j_{2}}
}    
\end{equation*}
be Cartesian. If
$j_{1}$ and $j_{2}$ are \'{e}tale
and a complex $\mathcal{U}\in D_{+}(X_{1})$,
then we have
an isomorphism
\begin{equation}\label{ebc}
j_{2}^{*}Rj_{1*}(\mathcal{U})
\simeq
Rj_{1*}^{\prime}j_{2}^{\prime*}(\mathcal{U})
\end{equation}
in $D_{+}(X_{2})$
by \cite[p.88, III, Theorem 1.15]{M}. 

Write $\Lambda=\mu_{p^{r}}^{\otimes n}$
and $U:=\mathfrak{X}\setminus Y$.
The natural inclusion
$j: U\hookrightarrow\mathfrak{X}$ is the complement of
the closed immersion
$\iota: Y\to \mathfrak{X}$.
For each $x\in U$,
$\bar{\{x\}}$ denotes its closure in $\mathfrak{X}$.
Then we have an isomorphism
\begin{equation*}
\HO^{s}_{x}(U_{\et}, \Lambda) 
\simeq
\HO^{s}_{\bar{\{x\}}}
(
\mathfrak{X}_{\et}, Rj_{*}\Lambda
)
\end{equation*}
for any integer $s$ by (\ref{ebc}). Thus the sequence
\begin{equation*}
\bigoplus_{x\in U^{(1)}}\HO_{x}^{n+1}(U_{\et}, \Lambda)
\to 
\HO^{n+1}_{\et}(U, \Lambda)
\to
\HO^{n+1}_{\et}(k(U
), \Lambda)
\to
\bigoplus_{x\in U^{(1)}}\HO_{x}^{n+2}(U_{\et}, \Lambda)
\end{equation*}
corresponds to the sequence
%
\begin{equation*}
\bigoplus_{x\in U^{(1)}}\HO_{\bar{\{x\}}}^{n+1}(\mathfrak{X}_{\et}, Rj_{*}\Lambda)
\to 
\HO^{n+1}_{\et}(\mathfrak{X}, Rj_{*}\Lambda)
\to
\HO^{n+1}_{\et}(k(\mathfrak{X}), \Lambda)
\to
\bigoplus_{x\in U^{(1)}}\HO_{\bar{\{x\}}}^{n+2}(\mathfrak{X}_{\et}, Rj_{*}\Lambda).
\end{equation*}
%
Put
\begin{math}
\mathcal{F}=\iota^{*}Rj_{*}\Lambda.
\end{math}

Consider the fiber product
\begin{equation*}
\xymatrix{
Y\ar[d]_-{\iota}
&\ar[l]_-{j_{3}^{\prime}}\ar[d]^-{\iota_{1}}
Y\cap\mathfrak{X}\setminus\bar{\{x\}}\\
\mathfrak{X}
&\ar[l]^-{j_{3}}
\mathfrak{X}\setminus\bar{\{x\}}
}    
\end{equation*}
where 
$j_{3}: \mathfrak{X}\setminus \bar{\{x\}}\to \mathfrak{X}$
is the open immersion of $\mathfrak{X}$.
Then we have an isomorphism
\begin{equation*}
(j_{3})^{*}\iota_{*}\mathcal{F}
\simeq
(\iota_{1})_{*}(j_{3}^{\prime})^{*}\mathcal{F}
\end{equation*}
by the proper base change theorem (\cite[pp.223--224, VI, Corollary 2.3]{M}) 
(or by \cite[p.69, II, Theorem 3.2]{M} and \cite[p.71, II, Corollary 3.5 (a)]{M}). 
So we have isomorphisms
\begin{align*}
\HO^{s}_{\et}(
\mathfrak{X}, 
R(j_{3})_{*}(j_{3})^{*}\iota_{*}\mathcal{F}
)
\simeq
\HO^{s}_{\et}(
Y,
R(j_{3}^{\prime})_{*}(j_{3}^{\prime})^{*}
\mathcal{F}
)
\end{align*}
and
\begin{align*}
\HO^{s}_{\bar{\{x\}}}(\mathfrak{X}_{\et}, \iota_{*}\mathcal{F})
\simeq
\HO^{s}_{\bar{\{x\}}\cap Y}(Y_{\et}, \mathcal{F})
\end{align*}
\normalsize
for any integer $s$. Hence the unit of the adjunction 
\begin{math}
\operatorname{id}
\to
\iota_{*}\iota^{*}
\end{math}
induces
a commutative diagram 
\begin{equation*}
\scriptsize{
\xymatrix
{
\displaystyle\bigoplus_{x\in U^{(1)}}
\HO^{n+1}_{\bar{\{x\}}}(\mathfrak{X}_{\et}, Rj_{*}\Lambda)
\ar[r]\ar[d]
&
\HO^{n+1}_{\et}(\mathfrak{X}, Rj_{*}\Lambda)
\ar[r]\ar[d]
&
\HO^{n+1}_{\et}(k(\mathfrak{X}),
\Lambda)
\ar[r]\ar[d]
&
\displaystyle\bigoplus_{x\in U^{(1)}}
\HO^{n+2}_{\bar{\{x\}}}(\mathfrak{X}_{\et}, Rj_{*}\Lambda)
\ar[d]
\\
\displaystyle
\bigoplus_{y\in Y^{(1)}}
\HO^{n+1}_{y}(Y_{\et}, \mathcal{F})
\ar[r]
&
\HO^{n+1}_{\et}(Y, \mathcal{F})
\ar[r]
&
\displaystyle
\bigoplus_{\eta\in Y^{(0)}}
\HO^{n+1}_{\et}(\kappa(\eta), \mathcal{F})
\ar[r]
&
\displaystyle
\bigoplus_{y\in Y^{(1)}}
\HO^{n+2}_{y}(Y_{\et}, \mathcal{F})
}
}
\end{equation*}
where the rows are exact.
Let $\tilde{\mathcal{O}_{\mathfrak{X}, \eta}}$ be the henselization of
the local ring 
$\mathcal{O}_{\mathfrak{X}, \eta}$ of $\eta\in Y^{(0)}$
and $K_{(\eta)}$ its fraction field.
Let $\tilde{j}: \spec(K_{(\eta)})\to \spec(\tilde{\mathcal{O}_{\mathfrak{X}, \eta}})$ 
be the complement of the closed immersion 
$\tilde{\iota}: \spec(\kappa(\eta))\to\spec(\tilde{\mathcal{O}_{\mathfrak{X}, \eta}})$.
Then we have an isomorphism
\begin{equation*}
\HO^{n+1}_{\et}(
\kappa(\eta), \mathcal{F}
)   
\simeq
\HO^{n+1}_{\et}(
\kappa(\eta), \tilde{\iota}^{*}
R\tilde{j}_{*}\Lambda
)
\end{equation*}
by (\ref{ebc}).
So we have isomorphisms
\begin{equation*}
\HO^{n+1}_{\et}(\kappa(\eta), \mathcal{F})
\simeq
\HO^{n+1}_{\et}(
\tilde{\mathcal{O}_{\mathfrak{X}, \eta}},
R\tilde{j}_{*}\Lambda
)
\simeq
\HO^{n+1}_{\et}(
K_{(\eta)},
\Lambda
)
\end{equation*}
by \cite[p.777, The proof of Proposition 2.2.b)]{Ge}.
Moreover,
we have an isomorphism
\begin{equation*}
\HO^{s-2}_{\et}(\kappa(x), \Lambda(-1))
\simeq
\HO^{s}_{x}(U_{\et}, \Lambda)
\end{equation*}
for any integer $s$
by the absolute purity theorem 
(\cite[p.241, VI, Theorem 5.1]{M}).
Hence we have a commutative diagram
with exact rows
\begin{equation}\label{bccom1}
\scriptsize{
\xymatrix
{
\displaystyle\bigoplus_{x\in U^{(1)}}
\HO^{n-1}_{\et}(\kappa(x), \Lambda(-1))
\ar[r]\ar[d]_{\gamma}
&
\HO^{n+1}_{\et}(\mathfrak{X}, Rj_{*}\Lambda)
\ar[r]\ar[d]_{\simeq}
&
\HO^{n+1}_{\et}(k(\mathfrak{X}), \Lambda)^{\prime}
\ar[r]\ar[d]
&
0
\\
\displaystyle
\bigoplus_{y\in Y^{(1)}}
\HO^{n+1}_{y}(Y_{\et}, \mathcal{F})
\ar[r]
&
\HO^{n+1}_{\et}(Y, \mathcal{F})
\ar[r]
&
\displaystyle
\bigoplus_{\eta\in Y^{(0)}}
\HO^{n+1}_{\et}(K_{(\eta)}, \Lambda)
&
}
}
\end{equation}
where 
\begin{equation*}
\HO^{n+1}_{\et}(k(\mathfrak{X}), \Lambda)^{\prime}    
=
\operatorname{Ker}\left(
\HO^{n+1}_{\et}(k(\mathfrak{X}), \Lambda)
\to
\displaystyle\bigoplus_{x\in U^{(1)}}
\HO^{n}_{\et}(
\kappa(x), 
\Lambda(-1))
\right).
\end{equation*}
The statement means that the right map in the diagram 
(\ref{bccom1}) is injective.
By the proper base change theorem, the middle map
in the diagram (\ref{bccom1})
is an isomorphism. 
So it suffices to show that the homomorphism 
\begin{equation*}
\gamma: 
\bigoplus_{x\in U^{(1)}}
\HO^{n-1}_{\et}(\kappa(x), \Lambda(-1))
\to
\bigoplus_{y\in Y^{(1)}}
\HO^{n+1}_{y}(Y_{\et}, \mathcal{F})
\end{equation*}
is surjective. 

Let $x\in U^{(1)}$. Since
\begin{math}
\bar{\{x\}}\to \mathfrak{X}\to
\spec (A)
\end{math}
is proper and quasi-finite, so is finite by 
\cite[p.6, I, Corollary 1.10]{M}. Hence
$\bar{\{x\}}$ is the spectrum of the henselian local domain 
\begin{math}
\mathcal{O}_{\mathfrak{X}, y}/\mathfrak{p}_{x}
\end{math}
by \cite[pp.32--33, I, Theorem 4.2]{M}
and \cite[p.34, I, Corollary 4.3]{M}.
Here $\mathfrak{p}_{x}$ is the prime ideal of $\mathcal{O}_{\mathfrak{X}, y}$
which corresponds to $x\in U^{(1)}$.
Therefore $\bar{\{x\}}$ meets $Y$ at one and only one point and 
\begin{equation*}
U^{(1)}=\bigsqcup_{y\in Y^{(1)}}
\left(
\spec 
(\mathcal{O}_{\mathfrak{X}, y})
\times_{\mathfrak{X}}
U
\right)^{(1)}.
\end{equation*}
Moreover we have isomorphisms
\begin{equation*}
\tilde{\mathcal{O}_{\mathfrak{X}, y}}/
\mathfrak{p}_{x}\cdot\tilde{\mathcal{O}_{\mathfrak{X}, y}}
\simeq
\widetilde{\mathcal{O}_{\mathfrak{X}, y}/\mathfrak{p}_{x}}
\simeq
\mathcal{O}_{\mathfrak{X}, y}/\mathfrak{p}_{x}
\end{equation*}
and
\begin{math}
\mathfrak{p}_{x}\cdot\tilde{\mathcal{O}_{\mathfrak{X}, y}}
\end{math}
is a prime ideal of $\tilde{\mathcal{O}_{\mathfrak{X}, y}}$. 
So we have a bijection
\begin{equation*}
P_{y}:=
(\spec(\tilde{\mathcal{O}_{\mathfrak{X}, y}})\times_{\mathfrak{X}}U)^{(1)} 
\to
(\spec(\mathcal{O}_{\mathfrak{X}, y})\times_{\mathfrak{X}}U)^{(1)}.
\end{equation*}
Hence the map $\gamma$ decomposes into a direct sum
\begin{equation*}
\gamma=\bigoplus_{y\in Y^{(1)}}
\left(
\gamma_{y}:
\bigoplus_{x\in P_{y}}\HO^{n-1}_{\et}(
\kappa(x), \Lambda(-1)
)
\to
\HO^{n+1}_{y}(Y_{\et}, \mathcal{F})
\right)
\end{equation*}
and it suffices to show that the map $\gamma_{y}$ 
is surjective. 

Put
\begin{equation*}
V_{y}=\spec(\tilde{\mathcal{O}_{\mathfrak{X}, y}})^{(1)}
\setminus(
\spec(\tilde{\mathcal{O}_{\mathfrak{X}, y}})\times_{\mathfrak{X}}U
)^{(1)}.
\end{equation*}
Let $(\tilde{\mathcal{O}_{\mathfrak{X}, y}})_{\eta}$
be the localization of
the henselian local ring
$\tilde{\mathcal{O}_{\mathfrak{X}, y}}$
at $\eta\in V_{y}$ and
$K_{(y, \eta)}$ the fraction field of
the henselization of
$(\tilde{\mathcal{O}_{\mathfrak{X}, y}})_{\eta}$.
Since
\begin{equation*}
\tilde{\mathcal{O}_{Y, y}}
=
\tilde{\mathcal{O}_{\mathfrak{X}, y}}
\otimes_{\mathcal{O}_{\mathfrak{X}, y}}\mathcal{O}_{Y, y}
\end{equation*}
by \cite[p.37, I, Examples 4.10 (c)]{M}, we have a commutative diagram with exact rows
\begin{equation}\label{grx}
\scriptsize{
\xymatrix
{
\HO^{n}_{\et}(K_{(y)}, \Lambda)
\ar[r]\ar[d]
&
\displaystyle\bigoplus_{x\in P_{y}}
\HO^{n-1}_{\et}(\kappa(x), \Lambda(-1))
\ar[r]\ar[d]_{\gamma_{y}}
&
\HO^{n+1}_{\et}(
\tilde{\mathcal{O}_{\mathfrak{X}, y}}\times_{\mathfrak{X}}U, \Lambda
)
\ar[r]\ar[d]_{\simeq}^{\varphi}
&
\HO^{n+1}_{\et}(K_{(y)}, \Lambda)^{\prime}
\ar[d]
\\
\displaystyle
\bigoplus_{\eta\in V_{y}}
\HO^{n}_{\et}(K_{(y, \eta)}, \Lambda)
\ar[r]
&
\HO^{n+1}_{y}(Y_{\et}, \mathcal{F})
\ar[r]
&
\HO^{n+1}_{\et}(\tilde{\mathcal{O}_{Y, y}}, \mathcal{F})
\ar[r]
&
\displaystyle
\bigoplus_{\eta\in V_{y}}
\HO^{n+1}_{\et}(K_{(y, \eta)}, \Lambda)
}
}
\end{equation}
where
$K_{(y)}$ is the fraction field of the henselian local ring 
$\tilde{\mathcal{O}_{\mathfrak{X}, y}}$
and
\begin{equation*}
\HO^{n+1}_{\et}(K_{(y)}, \Lambda)^{\prime}   
=
\operatorname{Ker}\left(
\HO^{n+1}_{\et}(K_{(y)}, \Lambda)
\to
\bigoplus_{x\in P_{y}}
\HO^{n}_{\et}(\kappa(x), \Lambda(-1))
\right).
\end{equation*}
In the diagram (\ref{grx}),
the left vertical map is surjective by Lemma \ref{WA}, 
the map $\varphi$ is an isomorphism by the proper base change theorem
and the right vertical map
is injective by Corollary \ref{LGLC}. Therefore the map $\gamma_{y}$
is surjective by chasing the diagram. 
This completes the proof.
\end{proof}
\begin{Remark}\upshape
In the proof of Theorem \ref{MTP}, we use 
the Bloch-Kato conjecture (\cite{V2}).
But it suffices to use \cite[p.131, Theorem (5.12)]{B-K} 
which is a special case
of the Bloch-Kato conjecture (cf.The proof of Lemma \ref{WA}).
\end{Remark}

\begin{thebibliography}{99}
\bibitem{SGA4}
\textsc{Artin, M., A. Grothendieck, and J.-L. Verdier.},
Th\'{e}orie des Topos et Cohomologie \'{E}tale des Sch\'{e}mas
(SGA 4), Tome 3.
Lecture Notes in Mathematics, vol. 305. Heidelberg, Berlin:
Springer-Verlag, 1973.
\bibitem{B}
\textsc{S.Bloch},
Algebraic cycles and Algebraic $K$-theory,
Adv. in Math. 61 (1986), no. 3, 267--304. 
\bibitem{B2}
\textsc{S.Bloch},
The moving lemma for higher Chow groups.
J. Algebraic Geom. \textbf{3} (1994), no.3, 537–-568.
\bibitem{B-K}
\textsc{S.Bloch and K.Kato},
$p$-adic \'{e}tale cohomology. 
Inst. Hautes \'{E}tudes Sci. Publ. Math. No. 63 (1986), 107--152.
\bibitem{B-O}
\textsc{S.Bloch and A.Ogus},
Gersten's conjecture and the homology of schemes. 
Ann. Sci. École Norm. Sup. (4) 7 (1974), 181--201 (1975). 
\bibitem{C}
\textsc{J.-L.Colliot-Th\'{e}l\`{e}ne},  
Quelques probl\`{e}mes locaux-globaux,
Personal Notes. (2011). 
\bibitem{C-H-K}
\textsc{
Colliot-Th\'{e}l\`{e}ne, Jean-Louis; 
Hoobler, Raymond T.; Kahn, Bruno},
The Bloch-Ogus-Gabber theorem,
Snaith, Victor P. (ed.), Algebraic K-theory. 
Papers from the 2nd Great Lakes conference, Canada, March 1996, 
in memory of Robert Wayne Thomason. 
Providence, RI: 
American Mathematical Society. Fields Inst. Commun. 16, 31--94 (1997).
\bibitem{G}
\textsc{K.Fujiwara},
A proof of the absolute purity conjecture (after Gabber), 
Algebraic geometry 2000, Azumino (Hotaka), 153--183, 
Adv. Stud. Pure Math., 36, Math.
Soc. Japan, Tokyo, 2002.
\bibitem{Ge}
\textsc{T.Geisser},
Motivic cohomology over Dedekind rings,
Math. Z. 248 (2004), no. 4, 773–-794. 
\bibitem{G-L}
\textsc{T.Geisser and M.Levine},
The $K$-Theory of fields in characteristic $p$,
Invent. Math. 139 (2000), no. 3, 459--493.
\bibitem{G-L2}
\textsc{T.Geisser and M.Levine},
The Bloch-Kato conjecture and a theorem of
Suslin-Voevodsky,
J. Reine Angew. Math.
530 (2001),
55--103.
\bibitem{G-S}
\textsc{M.Gros and N.Suwa},
La conjecture de Gersten pour les faisceaux de
Hodge-{W}itt logarithmique,
Duke Math. J. 57 (1988),
no. 2, 615--628.
\bibitem{H-H-K}
\textsc{D.Harbater and J.Hartmann and D.Krashen},
Local-global principles for {G}alois cohomology,
Comment. Math. Helv.89 (2014), 
no.1, 215–-253.
\bibitem{Hu}
\textsc{Y.Hu}, 
A Cohomological Hasse Principle Over Two-dimensional Local
Rings, Int. Math. Res. Not. IMRN 2017, no. 14, 4369--4397.
\bibitem{H}
\textsc{O.Hyodo},
A note on {$p$}-adic \'{e}tale cohomology in the semistable
reduction case,
Invent. Math. 91
(1988),
no.3,
543--557.
\bibitem{I}
\textsc{L.Illusie},
Complexe de de Rham-Witt et cohomologie cristalline,
Ann. Sci. \'{E}cole Norm. Sup. (4) 12 (1979),
no. 4, 501--661.
\bibitem{K}
\textsc{K.Kato},
A Hasse principle for two-dimensional global fields.
With an appendix by Jean-Louis Colliot-Th\'{e}l\`{e}ne.
J. Reine Angew. Math. \textbf{366} (1986), 142–-183.
\bibitem{L}
\textsc{M.Levine},
Techniques of localization in the theory of algebraic cycles,
J. Algebraic Geom. 10 (2001),
no.2,
299--363.
\bibitem{M}
\textsc{J.S.Milne},
Etale Cohomology,
Princeton Math. Ser. 
Princeton University Press, Princeton, N.J., 1980.
\bibitem{Q}
\textsc{D.Quillen},
Higher algebraic $K$-theory: I.
Higher K-theories (Proc.Conf., 
Battelle Memorial Inst., Seattle, Wash., 1972), pp. 85--147. 
Lecture Notes in Math.,
Vol. 341, Springer, 1973.
\bibitem{SM19}
\textsc{M.Sakagaito},
A note on Gersten’s conjecture for \'{e}tale
cohomology over two-dimensional henselian
regular local rings,
Comptes Rendus - Math\'{e}matique 358 (2020), no.1, 33--39.
\bibitem{SM}
\textsc{M.Sakagaito},
On a generalized Brauer group in mixed characteristic cases,
J. Math. Sci. Univ. Tokyo 27, 
No.1, 29--64 (2020). 
\bibitem{SaL}
\textsc{K.Sato},
Logarithmic Hodge--Witt sheaves on normal crossing
varieties,
Math. Z. 257 (2007), no. 4, 707--743.
\bibitem{SaP}
\textsc{K.Sato},
{$p$}-adic \'{e}tale {T}ate twists and arithmetic duality.
With an appendix by Kei Hagihara,
Ann. Sci. \'{E}cole Norm. Sup. (4)
40
(2007),
no.4,
519--588.
\bibitem{SaR}
\textsc{K.Sato},
Cycle classes for $p$-adic \'{e}tale 
Tate twists and the image 
of $p$-adic regulators.
Doc. Math. 18 (2013),
177--247.
\bibitem{V2}
\textsc{V.Voevodsky},
On motivic cohomology with $\mathbf{Z}/l$-coefficients,
Ann. of Math. (2) 174 (2011), no. 1, 401--438.
\end{thebibliography}
\end{document}